\documentclass{amsart}
\usepackage{amsmath}
\usepackage{amssymb}
\usepackage[all]{xy}
\usepackage{comment}
\usepackage{color}

\theoremstyle{plain}
\newtheorem{thm}{Theorem}[section]
\newtheorem{thm*}{Theorem}[section]
\newtheorem{cor}[thm]{Corollary}
\newtheorem{prop}[thm]{Proposition}
\newtheorem{lemma}[thm]{Lemma}

\newtheorem{lemma*}{Lemma}

\newtheorem{question}[thm]{Question}

\theoremstyle{definition}
\newtheorem{defn}[thm]{Definition}
\newtheorem{remark}[thm]{Remark}

\newtheorem{note}[thm]{Notation}

\newtheorem*{remark*}{Remark}

\newtheorem{ex}[thm]{Example}

\numberwithin{equation}{thm}

\def\Spec{\operatorname{Spec}\nolimits}

\def\Proj{\operatorname{Proj}\nolimits}

\newcommand{\bG}{\mathbb G}

\newcommand{\cO}{\mathcal O}

\newcommand{\bP}{\mathbb P}
\newcommand{\bZ}{\mathbb Z}

\newcommand{\cE}{\mathcal E}

\newcommand{\fp}{\mathfrak p}

\newcommand{\fg}{\mathfrak g}

\newcommand{\fq}{\mathfrak q}

\newcommand{\fm}{\mathfrak m}
\newcommand{\fn}{\mathfrak n}

\newcommand{\p}{\mathfrak p}

\newcommand{\ol}{\overline}
\newcommand{\ul}{\underline}

\def\Spec{\operatorname{Spec}\nolimits}
\def\sl2{\operatorname{SL_{2(2)}}\nolimits}
\def\Ga2{\operatorname{\mathbb G_{\rm a(2)}}\nolimits}

\newcommand{\Gr}{\mathbb G_{(r)}}

\newcommand{\bU}{\mathbb U}

\newcommand{\bu}{\bullet}

\date\today

\begin{document}

 \title[Support Theory for  Extended Drinfeld doubles]{Support Theory for Extended Drinfeld doubles}
 
 \author{Eric M. Friedlander}

\address {Department of Mathematics, University of Southern California,
Los Angeles, CA}
\email{ericmf@usc.edu}
-

\thanks{Partially supported by the Simons Foundation}



\begin{abstract} 
Following earlier work with Cris Negron on the cohomology of Drinfeld 
doubles $D(\mathbb G_{(r)})$, we develop a ``geometric theory" of support 
varieties for ``extended Drinfeld doubles" $\tilde D(\mathbb G_{(r)})$ 
of  Frobenius kernels $\mathbb G_{(r)}$ of smooth linear algebraic groups 
$\mathbb G$ over a field $k$ of characteristic $p > 0$.   To a        
$\tilde D(\mathbb G_{(r)})$-module $M$ we associate the space             
$\Pi(\tilde D(\mathbb G_{(r)}))_M$ of equivalence classes of ``pairs of 
$\pi$-points" and prove most of the desired properties of $M \mapsto 
\Pi(\tilde D(\mathbb G_{(r)}))_M$.   Namely, this association satisfies 
the ``tensor product property" and admits a natural continuous map 
$\Psi_{\tilde D}$  to cohomological support theory.  Moreover, for $M$ 
finite dimensional and with suitable conditions on $\mathbb G_{(r)}$, 
this association provides a ``projectivity test", $\Psi_{\tilde D}$ is 
a homeomorphism, and identifies $\Pi(\tilde D(\mathbb G_{(r)}))_M$ with 
the cohomological support variety of $M$ for various classes of 
$\tilde D(\mathbb G_{(r)})$-modules $M$.
\end{abstract}

\maketitle


\section{Introduction}

Cohomological support varieties for $kG$-modules for a finite group $G$ constitute natural extensions
of D. Quillen's geometric description of $\Spec H^\bu(G,k)$  \cite{Q1}, \cite{Q2}.   As introduced by J. Carlson
 in \cite{Ca}, the cohomological support variety $|G|_M$ of a finite dimensional $kG$-module $M$
 is the closed subvariety of  $\Spec H^\bu(G,k)$  associated to the (radical of) the kernel of the 
 natural map $H^\bu(G,k) \to Ext_G^*(M,M)$ of $k$-algebras.   To make this more geometric and more computable
 in the special case that $G$ is an elementary abelian $p$-group, J. Carlson conjectured and
 G. Avrunin and L. Scott proved (in \cite{A-S}) that $|E|_M$ can be identified with the ``rank variety"
 $V(E)_M$ of cyclic shifted subgroups along which the pull-backs of $M$ are not projective.  This
``geometric interpretation" of cohomological support varieties for modules for an elementary
abelian $p$-group was successively extended by 
B. Parshall and the author to modules for restricted Lie algebras \cite{F-Par}, by A. Suslin, C. Bendel, 
and the author in \cite{SFB2} to modules for infinitesimal group schemes, and to modules for arbitrary finite 
group schemes by J. Pevstova and the author \cite{F-P2}.  The
interweaving of the two approaches (one, primarily cohomological; the second, more geometric)
 has led to various constructions and invariants.
 
 	In this paper, we adapt the theory of $\pi$-points for finite group schemes to a theory of 
$\pi$-point pairs for a special class of finite dimensional Hopf algebras over $k$ which are neither
commutative nor cocommutative.  Namely, beginning with a smooth linear algebraic group $\bG$ over
$k$, we consider the Frobenius kernel $\bG_{(r)}$ (i.e., the kernel of the $r$-th iterate of the 
Frobenius map, $F^r: \bG \to \bG^{(r)}$ for some $r > 0$) and take the smash product algebra
$\tilde D(\Gr) \ \equiv \ k[\bG_{(r+1)}] \# k\Gr$ which maps surjectively to the Drinfeld double 
$D(\Gr)$ of $\Gr$.

	What enables our approach for $\tilde D(\Gr)$-modules
(and does not work for modules for the Drinfeld double $D(\Gr)$) is the existence of a 
sub Hopf algebra $O(\Gr) \subset \tilde D(\Gr)$ which essentially detects the cohomology of 
$\tilde D(\Gr)$ and which is isomorphic as an algebra to the ``group algebra" for a finite group
scheme.  Indeed, almost by definition, the $\pi$-point pair space $\Pi(\tilde D(\Gr))$ maps 
bijectively to the $\pi$-point pair space $\Pi(O(\Gr))$ which is homeomorphic to the $\pi$-point
space of the associated finite group scheme.  To establish the desired tensor product property
$\Pi(\tilde D(\Gr))_{M\otimes M^\prime} \ = \  \Pi(\tilde D(\Gr))_{M} \cap \Pi(\tilde D(\Gr))_{N}$
in Theorem \ref{thm:2coproducts},
we utilize this tensor product property for finite group schemes.  The key step requires the 
verification that this property remains valid upon changing coproduct structure on this group algebra
to the coproduct structure on $O(\Gr)$.  This invariance
of coproduct of supports of a tensor product is the content of Proposition \ref{prop:pullback}.
The tensor product property enables a natural topology on $\Pi(\tilde D(\Gr))$.

	In Section \ref{sec:cohom}, we extend the arguments of C. Negron and the author in \cite{F-N}
for the cohomology of Drinfeld doubles to the cohomology of the extended doubles $\tilde D(\Gr)$
which allows a comparison in Corollary \ref{cor:finite} of the cohomology of $\tilde D(\Gr)$ with that of 
the simpler subalgebra $O(\Gr)$.  In Proposition \ref{prop:cont}, we establish a commutative
square of continuous maps relating $\pi$-point pair spaces for $O(\Gr)$ and $\tilde D(\Gr)$ with
(projectivized) prime ideal spectra of their cohomology.

	In order to proceed further, we require the determination of the cohomology algebra 
$H^\bu(\tilde D(\Gr),k)$, the commutative subalgebra of the Ext-algebra $Ext_{\tilde D(\Gr)}^*(k,k)$ of classes
of even degree for $p > 2$ (and the full Ext-algebra if $p = 2$).  
This is achieved in Theorem \ref{thm:quasi-log} provided that 
$\bG$ admits a quasilogarithm (see Example \ref{ex:quasi}) and that $p^{r+1} > 2dim(\bG)$, 
again using techniques of \cite{F-N}.    This computation allows us to show that 
$\Psi_{\tilde D}: \Pi(\tilde D(\Gr)) \to \bP (H^\bu(\tilde D(\Gr),k))$ is a homeomorphism.  As a 
consequence, we conclude in Theorem \ref{thm:tD-test} that a finite dimensional $\tilde D(\Gr)$-module
is projective if and only if $\Pi(\tilde D(\Gr))_M$ is empty.

	The fact that $\Psi_{\tilde D}$ is a homeomorphism mapping $\Pi(\tilde D(\Gr))_M$ to 
$\bP (H^\bu(\tilde D(\Gr),k))_M$ does not a priori identify $\Pi(\tilde D(\Gr))_M$ with 
$\bP (H^\bu(\tilde D(\Gr),k))_M$.  In Section \ref{sec:relate}, we achieve this identification 
for various special classes
of $\tilde D(\Gr)$-modules, including Carlson's $L_\zeta$-modules, modules whose action
factors through the quotient $\tilde D(\Gr) \twoheadrightarrow k\Gr$ and modules whose 
restriction to $k[\bG_{(r+1)}]$ are projective.

	In Section \ref{sec:remarks}, we observe that the same techniques we use to  investigate 
$\tilde D(\Gr)$-modules apply to modules for closely related Hopf algebras.  We also mention 
a couple of questions worthy of further investigation.

	We thank both Cris Negron and Julia Pevtsova for their contributions to the 
foundations of this work, and the Institute for Advanced Study for its (virtual) hospitality.

\vskip .2in


\section{ $ D(\Gr), \ \tilde D(\Gr), \ O(\Gr)$,  $ \Gr \times$ (\protect\ul{$\fg$}$^{(r)})_{(1)} $}


Throughout this work, $\bG$ will denote a connected linear algebraic group smooth over $k$. 
If $A$ is a $k$-algebra and $K/k$ is a field extension, 
then we denote by $A_K$ the base change $K\otimes A$;
for an $A$-module $M$, we use the notation $M_K$ for the $A_K$-module $K \otimes M$.
For a finite dimensional $k$-algebra $A$, we denote by $H^\bu(A,k) \ \subset Ext_A^*(k,k)$
the commutative subalgebra of even dimensional classes if $p > 2$ and the full algebra
$Ext_A^*(k,k)$ if $p=2$.

For any $r > 0$, we denote by $F^r: \bG \to \bG^{(r)}$ the $r$-th iterate of the
Frobenius map $F: \bG \to \bG^{(1)}$ given by the $k$-linear map 
$$F^*: k[\bG^{(1)}] \equiv k\otimes_\phi k[\bG] \to k[\bG] , \quad \lambda \otimes f \mapsto \lambda f^p$$ 
of coordinate algebras (where $k\otimes_\phi -$ denotes the base change along the $p$-th power map $\phi:k \to k$).
The {\it $r$-th Frobenius kernel} of $\bG$ is the finite group scheme 
\begin{equation}
\label{eqn:Gr}
\Gr \ \equiv \ ker\{ F^r: \bG \to \bG^{(r)} \} \ \hookrightarrow \bG,
\end{equation} 
whose coordinate algebra $k[\Gr]$ is the local commutative algebra defined as the quotient of the localization
$(k[\bG])_{(e)}$ of $k[\bG]$ at the identity $e \in \bG$ modulo the ideal generated by elements $x^{p^r}, x \in \fm_e$.
For example, 
\ $k[GL_{n(r)}] \ = \ k[x_{i,j}]/(x_{i,j}^{p^r}-\delta_{i,j})$.  

The coordinate algebra $k[\bG]$ of $\bG$ is a Hopf algebra with comultiplication determined by the 
multiplication of $\bG$.   For a finite group scheme $G$ such as $\bG_{(r)}$, we consider the $k$-linear dual
$kG$ of $k[G]$ which is a finite dimensional, co-commutative Hopf algebra (over $k$); we refer to $kG$ as the
group algebra of $G$.   For such a finite group scheme $G$, a (rational) $G$-module is a 
$k$-vector space equipped with a $k$-linear module structure, $kG \times M \to M$. 
We shall frequently refer to this structure  as the ``action of $G$ on $M$".

We recall that the Drinfeld double of a Hopf algebra $H$ is the double crossed product 
of $H^{*cop} $ and $H$, where $H^{*cop}$ has the same algebra structure as $H^*$ 
and coproduct given by the opposite of the coproduct of $H^*$.
If $H$ is finite dimensional and cocommutative, then $D(H)$ is isomorphic to the smash product $H^{*op} \# H$
(see \cite[10.3.10]{Mont}).

In Hopf algebra standard notation, the coadjoint action of an element $h \in H$ on an element $f \in H^*$ is given
by $\sum f(h_1 \rightharpoonup f \leftharpoonup s^{-1}(h_2) )$ where $s$ is the antipode of $H$, $\Delta$ is the coproduct
of $H^*$ and $\Delta(h)$ is written as $\sum h_1 \otimes h_2$.  By a result of D. Radford \cite[10.3.5]{Rad}, 
the product structure of $D(H)$ is given for finite dimensional $H$ by 
$$(f\#h)(f^\prime\#h^\prime) \ = \ \sum f(h_1 \rightharpoonup f^\prime \leftharpoonup s^{-1}(h_2) )\# h_3 h^\prime.$$
If $H$ is both finite dimensional and cocommutative, this multiplication is that of the smash product algebra
$H^* \# H$:
$$(f\#h)(f^\prime\#h^\prime) \ = \ \sum f h_1(f)\# h_2 h^\prime;$$
the coproduct of $D(H)$ is the tensor product of the opposite coproduct of $H^*$ and the coproduct of $H$.

\begin{defn}
\label{defn:DD}
As discussed above, the Drinfeld double of  the finite dimensional cocommutative Hopf algebra $k\Gr$ 
has algebra structure given by the smash product
\begin{equation}
\label{eqn:Drinfeld}
D(\bG_{(r)}) \ \equiv \ k[\bG_{(r)}]  \# k\bG_{(r)}.
\end{equation}
with respect to the 
right coadjoint action of $\bG_{(r)}$ on $k[\bG_{(r)}]$ and coalgebra structure given by the tensor product of
the co-opposite cooproduct on $k[\Gr]$ (i.e., the opposite product structure of $k\Gr$) and the coproduct of $k\Gr$.

The extended Drinfeld double $\tilde D(\Gr) $ is the analogous smash product with respect to the 
restriction to $\Gr$ of the coadjoint action  of $\bG_{(r+1)}$ on $k[\bG_{(r+1)}]$:
\begin{equation}
\label{eqn:extDrinfeld}
\tilde D(\bG_{(r)}) \ \equiv \ k[\bG_{(r+1)}]  \# k\bG_{(r)} \  \hookrightarrow \ D(\bG_{(r+1)}).
\end{equation}
\end{defn}

\vskip .1in
To verify that (\ref{eqn:extDrinfeld}) is a well defined inclusion of Hopf algebras, we observe that the restriction of
the coadjoint action of $\bG_{(r+1)}$ on $k[\bG_{(r+1)}]$ along $\Gr \hookrightarrow \bG_{(r+1)}$ is
the action used to define $\tilde D(\Gr)$, then check directly that the multiplication and comultiplication of
$D(\bG_{(r+1)})$ restricts via the evident inclusion of $k$-vector spaces
$k[\bG_{(r+1)}] \otimes k\Gr \ \subset \ k[\bG_{(r+1)}] \otimes k\bG_{(r+1)}$ to the  multiplication 
and comultiplication of $\tilde D(\Gr)$.

\vskip .1in

We shall be interested in (left) modules for the extended Drinfeld double $\tilde D(\Gr)$ which can be 
viewed as $k$-vector spaces $N$ with an action $\tilde D(\Gr) \times N \to N$.  The condition on such a  
pairing to be a left module structure can be expressed as the restriction to elements in 
$\tilde D(\Gr)  \subset D(\bG_{(r+1)})$ of the condition that a pairing $D(\bG_{(r+1)}) \times N \to N$ be 
a $D(\bG_{(r+1)})$ -module structure.  This latter condition is made explicit in the following proposition
(with $H = k\bG_{(r+1)}$).

We first observe that for any Hopf algebra $H$ there are natural algebra embeddings of $H, H^*$ in $D(H)$: 
$h \in H \mapsto 1\#h; \ f \in H^* \mapsto f\#1$.

\begin{prop} \cite[10.6.16]{Mont}, \cite{Maj}
Let $H$ be a finite dimensional Hopf Algebra.
The data of a $D(H)$-module structure on a vector space $N$ consists of module structures on $N$
for both $H^*$ and $H$ such that 
$$h \cdot (f\cdot n) \ = \ \sum (h_1 \rightharpoonup f_2) \cdot ((h_2 \leftharpoonup f_1)\cdot n) \ \in \ N.$$

This is equivalent to a left/right Yetter-Drinfeld module structure on $N$ which consists of the data of a $k\Gr$-module
structure and a right $k\Gr$-comodule structure on $N$ satisfying the Yetter-Drinfeld condition:
$$ \sum h_1 \cdot n_o \otimes h_2n_1 \ = \ \sum(h_1n)_0 \otimes (h_2\cdot n)_1h_1 \ \in \ N \otimes k\Gr.$$
Here, we have written the coaction $N \to N\otimes H$ as $n \mapsto \sum n_0 \otimes h_1$.
\end{prop}

\vskip .1in

We relate the extended Drinfeld double $\tilde D(\Gr)$ to various other structures.

\begin{defn}
\label{defn:OGr}
We define 
\begin{equation}
\label{eqn:OGr}
O(\Gr) \ \equiv \ k[(\bG^{(r)})_{(1)}] \otimes k\Gr.
\end{equation}
Using the natural identification of $(\bG^{(r)})_{(1)}$ with $\bG_{(r+1)}/\Gr$, we obtain the natural 
embedding
\begin{equation}
\label{eqn:O}
i_O: O(\bG_{(r)})  \ \hookrightarrow \ \tilde D(\bG_{(r)}).
\end{equation}
 \end{defn}
 
\vskip .1in

\begin{prop}
\label{prop:OG}
The embedding $\Gr \ \subset \bG_{(r+1)}$ determines a $\Gr$-equivariant quotient map 
$k[\bG_{(r+1)}] \twoheadrightarrow k[\Gr]$ of Hopf algebras and thus a quotient map of Hopf algebras
\begin{equation}
\label{eqn:quot}
q: \tilde D(\Gr) \ \twoheadrightarrow \ D(\Gr).
\end{equation}

Similarly, the quotient map $\bG_{(r+1)} \twoheadrightarrow (\bG^{(r)})_{(1)}$ determining $i_O$ of 
(\ref{eqn:O}) is an embedding of Hopf algebras

The coproduct structure inherited by $O(\Gr)$ is the tensor product of the co-opposite of the coproduct on $k[\bG^{(r)})_{(1)}]$
(in other words the dual of the opposite multiplication of $k(\bG^{(r)})_{(1)}$) and the coproduct of $ k\bG_{(r)}$
(the dual of the multiplication of $k[\Gr]$).  Thus, $O(\Gr)$ is co-commutative if and only if $\fg = Lie(\bG)$ is commutative.
\end{prop}

\begin{proof}
Since each $\bG_{(r)}$ is normal in $\bG$, the embedding $\Gr \ \subset \bG_{(r+1)}$ is $\bG$-equivariant
so that $k[\bG_{(r+1)}] \twoheadrightarrow k[\bG_{(r+1)}/\bG_{(r)}]$ is a $\bG$-equivariant quotient map of Hopf algebras;
thus, this $\Gr$-equivariant quotient determines the map of Hopf algebras $q: \tilde D(\Gr) \ \twoheadrightarrow \ D(\Gr)$.

Similarly, the quotient $\bG_{(r+1)} \twoheadrightarrow \bG^{(r)})_{(1)}$ is also $\bG$-equivariant,
determining the $\bG$-equivariant embedding $k[\bG^{(r)})_{(1)}] \hookrightarrow k[\bG_{(r+1)}]$ of Hopf algebras.
Observe that $\Gr \subset \bG$ acts trivially on $k[\bG^{(r)})_{(1)}]$, so that we obtain 
the embedding (as algebras) of the tensor product \ $k[\bG^{(r)})_{(1)}]\otimes k\Gr  \  \hookrightarrow \ \tilde D(\Gr)$.
\end{proof}

\begin{prop}
\label{prop:Oflat}
The subalgebra $k[\bG^{(r)})_{(1)}] \ \simeq \ k[\bG_{(r+1)}/\Gr]$ of the commutative algebra $k[\bG_{(r+1)}]$ gives $k[\bG_{(r+1)}]$
the structure of a free $k[\bG^{(r)})_{(1)}]$-module with basis $\{ \tilde f_i \}$ given by lifting along 
$k[\bG_{(r+1)}] \to k[\bG_{(r)}]$ a $k$-basis $\{ f_i \}$ for $k[\bG_{(r)}]$. 
 
This structure determines a free (right) $O(\Gr)$-module structure on $\tilde D(\Gr)$ with basis $\{ \tilde f_i \# 1\}$.
\end{prop}

\begin{proof}
We choose a $k$-linear section $\sigma: k[\Gr] \to k[\bG_{(r+1)}]$ of the quotient map $k[\bG_{(r+1)}] \twoheadrightarrow k[\Gr]$
and consider the $k[\bG_{(r+1)}/\Gr]$-bilinear map
\begin{equation}
\label{eqn:sig}
k[\Gr] \otimes k[\bG_{(r+1)}/\Gr] \ \to \ k[\bG_{(r+1)}], \quad f_i \otimes h \mapsto \sigma(f_i)\cdot h.
\end{equation}
We readily check that this induces a $k$-linear isomorphism, for example by showing that tensors 
of elements of $k$-bases for $k[\Gr]$ and $ k[\bG_{(r+1)}/\Gr]$ determine 
a $k$-linearly independent set in $k[\bG_{(r+1)}]$.

Using the fact that $k[\bG_{(r+1)}/\Gr]\otimes 1 \ \subset \tilde D(\Gr)$ is central, we readily check that
$$k[\Gr] \otimes O(\Gr) \to \ \tilde D(\Gr), \quad f \otimes (f^\prime \otimes h) \mapsto 
(\sigma(f)\#1) \cdot (f^\prime \# h)$$
is an $O(\Gr)$-linear isomorphism
\end{proof}

The following observation about $O(\Gr)$ will prove central in our investigation of $\tilde D(\Gr)$-modules.
The reader should not confuse the finite group scheme  $ (\ul \fg^{(r)})_{(1)}$ occurring in the proposition below 
with $(\bG^{(r)})_{(1)}$.

\begin{prop}
\label{prop:ulg}
Let $\ul \fg$ denote the vector group scheme $\Spec S^\bu(\fg^*)$ (ignoring the Lie algebra structure on $\fg$). 
Then the group algebra of the finite group scheme $(\ul \fg^{(r)})_{(1)} \times \bG_{(r)}$ is isomorphic
as an algebra  to $\cO(\Gr)$.

On the other hand, the dual of the Hopf algebra $O(\Gr)$ is not commutative (unlike the coordinate algebra
of  $ (\ul \fg^{(r)})_{(1)} \times \bG_{(r)}$).
\end{prop}

\begin{proof}
Essentially by definition, elements $X^{(r)} \in Lie(\bG^{(r)})$ are distributions at the identity of $\bG^{(r)}$ 
of $k[\bG^{(r)}]$ of order 1 (see \cite[I.7.10]{J}),  This identification determines
an embedding of $\fg^{(r)*}$ into the maximal ideal of $k[\bG_{(r+1)}/\bG_{(r)}]$ which necessarily sends 
each $(X^{(r)})^p \in S^p(\fg^{(r)*})$ to 0.  The resulting map 
$$\nu: S^\bu(\fg^{(r)*})/((X^{(r)})^p,X\in \fg^*) \ \to \ k[\bG_{(r+1)}/\bG_{(r)}]$$
is injective, and thus surjective by dimension reasons. 
The proof is completed by observing that the group algebra of  $(\ul \fg^{(r)})_{(1)}$ equals
$S^\bu(\fg^{(r)*})/((X^{(r)})^p,X\in \fg^*)$.
\end{proof}

We conclude this section with some remarks about  $\pi$-point spaces $\Pi(G_1 \times G_2)$ of
a product of finite group schemes which we shall use for $G_1 \times G_2 = (\ul \fg^{(r)})_{(1)} \times \bG_{(r)}$.
Recall that a $\pi$-point $\alpha_K: K[t]/t^p \to kG$ of a finite group scheme $G$ over $k$ is a 
flat map of algebras which factors through the group algebra of some unipotent abelian subgroup scheme
$C_K \subset G_K$, where $K$ is some field extension of $k$.  Thus, a $\pi$-point of $G_1 \times G_2$ 
is a pair of $\pi$-points $\alpha_{K,1}: K[t]/t^p \to KG_1, \ \alpha_{K,2}: K[t]/t^p \to KG_2$ corresponding
to a flat map 
$$\alpha_{K,1}+\alpha_2: K[t]/t^p \to KC_1 \otimes KC_2, \quad t \mapsto \alpha_{K,1}(t)\otimes 1 + 1 \otimes 
\alpha_{K,1}(t)$$
with $C_{K,i} \subset G_{K,i}$ unipotent abelian subgroup schemes for $i= 1,2$.  

	It is important here for $\alpha_{K,1}(t), \ \alpha_{K,1}(t)$ to commute in order that $\alpha_{K,1}+\alpha_{K,2}$ is well defined.

By \cite[Thm 3.6]{F-P2}, two $\pi$-points $\alpha_K: K[t]/t^p \to KG, \beta_L: L[t]/t^p \to LG$ 
are equivalent if and only if the two ideals of $H^\bu(G,k)$
$$ker\{ \alpha_K^*: H^\bu(G,K) \to H^*(K[t]/t^p,K)\} \ \cap \ H^\bu(G,k),$$
$$ker\{ \beta_L^*: H^\bu(G,L) \to H^*(L[t]/t^p,L)\} \ \cap \ H^\bu(G,k)$$
are equal.  In other words, this gives a natural bijection of sets (which is sharpened to an isomorphism of
schemes in \cite{F-P2})
\begin{equation}
\label{eqn:PiG}
\Pi_G: \Pi(G)  \ \stackrel{\sim}{\to} \ \bP H^\bu(G,k).
\end{equation}

Observe that the flat map $(\alpha_{K,1},\alpha_{K,2}): K[t]/t^p \to KG_1 \times KG_2$
induces the map 
$$(\alpha_{K,1}^* \otimes 1) + (1 \otimes \alpha_{K,2})^*: H^*(C_1,K) \otimes_K H^\bu(C_2,K) \to H^*(K[t]/t^p,K).$$
We denote the composition of this map with the restriction map
$$(\alpha_{K,1 }+ \alpha_{K,2})^*: H^*(G_1 \times G_2,K) \ \to \ H^*(K[t]/t^p,K).$$

Thus, we have proved the following proposition.

\begin{prop}
\label{prop:equal-prod}
Two $\pi$-points $(\alpha_{K,1},\alpha_{K,2})$ and $(\beta_{L,1},\beta_{L,2})$ of the finite group 
scheme $G_1 \times G_2$ are equivalent if and only if the two ideals of $H^\bu(G_1\times G_2,k)$
$$ker\{ (\alpha_{K,1 }+ \alpha_{K,2})^*: H^\bu(G_1\times G_2,K) \to H^*(K[t]/t^p,K)\} \ \cap \ H^\bu(G_1 \times G_2,k),$$
$$ker\{(\beta_{L,1 }+ \beta_{L,2})^*: H^\bu(G_1 \times G_2,L) \to H^*(L[t]/t^p,L)\} \  \cap \ H^\bu(G_1 \times G_2,k)$$
are equal. 

In particular, if $\alpha_{K_i}$ is equivalent to $\beta_{L,i}$ for $i= 1,2$, then 
$(\alpha_{K,1},\alpha_{K,2})$ and $(\beta_{L,1},\beta_{L,2})$ are equivalent.
\end{prop}

\vskip .2in


\section{The space $ \Pi(\tilde D(\Gr))$ of $\pi$-point pairs}
\label{sec:pairs}

In this section we introduce an extension of ``rank varieties" which applies to $\tilde D(\Gr)$-modules.
Our construction is a modification of the ``$\pi$-point construction" of \cite{F-P2}.   The strategy we 
follow is relatively straight-forward.   In view of the homological computations of \cite{F-N}, it 
seems natural to use pairs of $\pi$-points $(\alpha_K,\beta_K)$ with $\alpha_K: K[t]/t^p \to K[\bG_{r+1)}]$
and $\beta_K: K[T]/t^p \to K\Gr$.  Yet the most natural way to do this (which would also apply to 
$D(\Gr)$) appears to fail because $\alpha_K(t), \beta_K(t)$ need not commute.  We work around this
difficulty by restricting $\tilde D(\Gr)$-modules to $O(\Gr)$.

\vskip .1in
\begin{defn}
\label{defn:pairs}
A $\pi$-point pair $(\alpha_K,\beta_K)$ of $\tilde D(\bG_{(r)})$
is a pair of maps of $K$-algebras for some field extension $K/k$,
$$\alpha_K: K[t]/t^p \to K[\bG_{(r+1)}/\bG_{(r)}] \hookrightarrow 
K[\bG_{(r+1)}], \quad \beta_K: K[t]/t^p \to K\bG_{(r)},$$
such that $\beta_K$ either sends $t$ to 0 or is a $\pi$-point of $\bG_{(r)}$, $\alpha_K$ either 
sends $t$ to 0 or is flat,
and at least one of $\alpha_K, \beta_K$ is flat.
\end{defn}

\vskip .1in

The following lemma provides an essential property of a $\pi$-point pair $(\alpha_K,\beta_K)$. 

\begin{lemma}
\label{lem:flat}
A  $\pi$-point pair  $(\alpha_K,\beta_K)$ of $\tilde D(\bG_{(r)})$ determines a flat map
$$\alpha_K + \beta_K: K[t]/t^p \ \to \ O(\Gr)_K, \quad t \mapsto \alpha_K(t) \# 1 + 1\# \beta_K(t)$$
whose composition with $i_O: O(\bG_{(r),K}) \  \hookrightarrow  \ \tilde D(\bG_{(r),K})$ is also flat.
\end{lemma}

\begin{proof}
We consider the composition 
$$K[t] \ \stackrel{\Delta}{\to} \ K[t]/p^{\otimes 2} \ \stackrel{\alpha_K \boxtimes\beta_K}{\to} 
\ K[\bG_{(r+1)}/\Gr] \otimes K\Gr = O(\Gr)_K,
\quad t \mapsto \alpha_K(t) \otimes \beta_K(t).$$
Since $\alpha_K(t),\beta_K(t)$ commute in $O(\Gr)_K$, this composition factors through $\alpha_K + \beta_K: K[t]/t^p \ \to \ O(\Gr)_K.$
Since $\Delta: K[t]/t^p \ \stackrel{\Delta}{\to} \ K[t]/p\otimes K[t]/t^p$ is flat, since both $\alpha_K$ and $\beta_K$ are either flat or ``trivial", since $O(\Gr)$ is flat over each of its tensor factors,
and since $\alpha_K \boxtimes\beta_K$ is flat if $\alpha_K(t) \not= 0 \not= \beta_K(t)$, 
we conclude that  $\alpha_K + \beta_K: K[t]/t^p \ \to \ O(\Gr)$
is flat.  By Proposition \ref{prop:Oflat}, this implies the composition $i_O \circ (\alpha_K + \beta_K): K[t]/t^p \ \to \ \tilde D(\Gr)$ is also flat.
\end{proof}

The following definition of equivalence for $\pi$-point pairs \ $(\alpha_K,\beta_K), (\alpha^\prime_L,\beta^\prime_L)$ \
is motivated by Proposition \ref{prop:equal-prod}.

\begin{defn}
\label{defn:equiv}
Two $\pi$-point pairs  $(\alpha_K,\beta_K), (\alpha^\prime_L,\beta^\prime_L)$ of $\tilde D(\Gr)$ 
are said to be equivalent (respectively, $O(\Gr)$-equivalent) if there exists some common field extension 
$\Omega$ of both $K$ and $L$  such that  for all finite dimensional $\tilde D(\bG_{(r)})$-modules $M$ 
(resp., all finite dimensional $O(\bG_{(r)})$-modules $M$)
$(\alpha_\Omega + \beta_\Omega)^*(M_\Omega)$ is projective as an $\Omega[t]/t^p$-module if and only if 
$(\alpha_\Omega^\prime + \beta_\Omega^\prime)^*(M_\Omega)$ is projective.

We denote by $\Pi (\tilde D(\bG_{(r)}))$ the set of  equivalence classes of $\pi$-point 
pairs of $\tilde D(\bG_{(r)})$, and we denote by $\Pi (O(\bG_{(r)}))$ the set of  
$O(\Gr)$-equivalence classes of $\pi$-point pairs of $\tilde D(\bG_{(r)})$.
For any $\tilde D(\Gr)$-module $M$, we define the {\it $\Pi$-pair set}
$$\Pi (\tilde D(\bG_{(r)}))_M \quad \subset \quad \Pi (\tilde D(\bG_{(r)}))$$
to be the subset of equivalence classes  of $\pi$-point pairs of $\tilde D(\bG_{(r)})$ such that $(\alpha_K+\beta_K)^*(M)$
is not projective.   Similarly, for any $O(\Gr)$-module $M^\prime$, we define \\
$\Pi (O(\bG_{(r)}))_{M^\prime} \subset \Pi (O(\bG_{(r)}))$
to be the subset of equivalence classes  of $\pi$-point pairs of $\tilde D(\bG_{(r)})$ such that $(\alpha_K+\beta_K)^*(M^\prime)$
is not projective.
\end{defn}

\vskip .1in

%

\begin{prop}
\label{prop:Pi-equiv}
The natural surjection $\Pi_{\tilde D}: \Pi (O(\bG_{(r)})) \ \twoheadrightarrow \ \Pi (\tilde D(\bG_{(r)})) $ is a bijection.

Moreover, for any finite dimensional $\tilde D(\Gr)$-module $M$, $\Pi_{\tilde D}$ restricts to a bijection
$$\Pi_{\tilde D,M}:  \Pi (O(\bG_{(r)}))_M  \ \stackrel{\sim}{\to} \ \Pi (\tilde D(\bG_{(r)}))_M.$$
\end{prop}

\begin{proof}
If $(\alpha_K,\beta_K), (\alpha^\prime_L,\beta^\prime_L)$ are $O(\bG_{(r)})$-equivalent, then they are 
necessarily equivalent since equivalence involves the same condition as $O(\bG_{(r)})$-equivalence except 
that it only requires consideration of finite dimensional $\tilde D(\bG_{(r)})$-modules.  Thus, we have well defined, 
surjective map $ \Pi (O(\bG_{(r)})) \ \twoheadrightarrow \ \Pi (\tilde D(\bG_{(r)}))$.

To prove injectivity, let $(\alpha_K,\beta_K), (\alpha^\prime_L,\beta^\prime_L)$ represent distinct points
in $\Pi (\tilde D(\bG_{(r)}))$ and choose some finite dimensional $O(\Gr)$-module $N$ such that 
$(\alpha_\Omega + \beta_\Omega)^*(N_\Omega)$ is projective whereas 
$(\alpha_\Omega^\prime + \beta_\Omega^\prime)^*(N_\Omega)$
is not projective for some some common field extension $\Omega$ of both $K$ and $L$.   Consider the
$\tilde D(\Gr)$-module $M \ \equiv \ \tilde D(\Gr) \otimes_{O(\Gr)} N$ and set $N^\prime \equiv 
 i_O^*(M)$.  By Proposition \ref{prop:Oflat}, $N^\prime$ 
 is isomorphic to a direct sum of copies of $N$ (indexed by a $k$-basis of $k[\Gr]$).
 Thus, $(\alpha_\Omega + \beta_\Omega)^*(N^\prime_\Omega)$ is projective whereas 
$(\alpha_\Omega^\prime + \beta_\Omega^\prime)^*(N^\prime_\Omega)$ is not projective.
Since $N^\prime$ is the restriction of the $D(\Gr)$-module $M$, we conclude that the images of 
$(\alpha_K,\beta_K), (\alpha^\prime_L,\beta^\prime_L)$ in $\Pi(\tilde D(\Gr))$ are distinct.

Finally, the bijectivity of $\Pi_{\tilde D,M}$
follows immediately from the definitions and the bijectivity of $\Pi_{\tilde D}$.
\end{proof}

We next observe that the condition of equivalence of $\pi$-point pairs is equivalent to the 
seemingly stronger equivalence relation obtained by dropping the finite dimensionality condition
on $M$ in Definition \ref{defn:equiv}.

\begin{prop}
\label{prop:stronger}
Two $\pi$-point pairs $(\alpha_K,\beta_K), (\alpha^\prime_L,\beta^\prime_L)$ of $\tilde D(\Gr)$ 
(respectively $O(\Gr)$) are equivalent if and only if there exists some common field extension 
$\Omega$ of both $K$ and $L$  such that  for all $\tilde D(\bG_{(r)})$-modules $M$ 
(resp., all $O(\bG_{(r)})$-modules $M$)
$(\alpha_\Omega + \beta_\Omega)^*(M_\Omega)$ is projective as an $\Omega[t]/t^p$-module if and only if 
$(\alpha_\Omega^\prime + \beta_\Omega^\prime)^*(M_\Omega)$ is projective.
 \end{prop}
 
 \begin{proof}
Denote by $\tilde \Pi(\tilde D(\Gr))$ (respectively, $\tilde \Pi(O(\Gr))$) the set of equivalence classes
of $\pi$-point pairs of $\tilde D(\Gr)$ (resp., $O(\Gr)$)  using the equivalence relations of the statement
of this proposition.  We readily verify that we have a commutative square of sets
\begin{equation}
\begin{xy}*!C\xybox{%
\xymatrix{
\tilde \Pi(O(\Gr))\ar[r]^-{\Pi_{\tilde D}} \ar[d]_{\tilde \Pi_O} & \tilde \Pi(\tilde D(\Gr)) \ar[d]^{\tilde \Pi_D} \\
\Pi(O(\Gr)) \ar[r]^-{\Pi_D} & \Pi(\tilde D(\Gr)).}
}\end{xy}
\end{equation}
The map $\Pi_D$ is a bijection as seen in Proposition \ref{prop:Pi-equiv}; the proof of that
proposition also proves that $\Pi_{\tilde D}$ is a bijection.    The fact that $\tilde \Pi_O$
is a bijection follows from  \cite[Thm 4.6]{F-P2} applied to the finite group scheme 
$(\ul \fg^{(r)})_{(1)} \times \bG_{(r)}$ (see Proposition \ref{prop:OtoG} below).  
Thus, the fact that $\tilde \Pi_D$ is a bijection
follows from the surjectivity of $\Pi_{\tilde D}$ and a simple diagram chase.
 \end{proof}

The following properties of $M \ \mapsto \ \Pi (\tilde D(\bG_{(r)}))_M$ are each verified 
by restricting along flat maps $(\alpha_\Omega + \beta_\Omega)^*$ associated to 
$\pi$-point pairs  $(\alpha_K,\beta_K)$, thereby reducing each assertion to the 
corresponding assertion about $K[t]/t^p$-modules.

\begin{prop}
\label{prop:basic}
Let $M_1, \ M_2, \ M_3$ be $\tilde D(\Gr)$-modules.  
\begin{enumerate}
\item
$\Pi (\tilde D(\bG_{(r)}))_k \ = \ \Pi (\tilde D(\bG_{(r)}))$.
\item
If $M_1$ is a projective $\tilde D(\Gr)$-module, then $\Pi (\tilde D(\bG_{(r)}))_{M_1} \ = \ \emptyset$.
\item
If $0 \to M_1 \to M \to M_2 \to$ is exact and if $\sigma$ is a permutation of $\{1,2,3 \}$, then 
$$\Pi (\tilde D(\bG_{(r)}))_{M_{\sigma(1)}} \ \subset \ \Pi (\tilde D(\bG_{(r)}))_{M_{\sigma(2)}} \cup \Pi (\tilde D(\bG_{(r)}))_{M_{\sigma(3)}}.$$ 
\item
$$\Pi (\tilde D(\bG_{(r)}))_{M_1 \oplus M_3} \ = \ \Pi (\tilde D(\bG_{(r)}))_{M_1} \cup \Pi (\tilde D(\bG_{(r)}))_{M_3}.$$
\end{enumerate}
\end{prop}

\vskip .1in

We proceed to show that $\Pi(O(\Gr))$ can be ``identified" with the set of Zariski points of the $\pi$-point scheme 
of the finite group scheme $(\ul \fg^{(r)})_{(1)}\times \bG_{(r)}$ introduced below.

\begin{prop}
\label{prop:OtoG}
Let $\ul{\fg} \ \simeq \ (\bG_a)^{\times dim(\fg)}$ denote the vector group scheme associated to the underlying
vector space of $\fg = Lie(\bG)$.  As an algebra $O(\Gr)$ is isomorphic to the group algebra of the 
infinitesimal group scheme $(\ul \fg^{(r)})_{(1)}\times \bG_{(r)}$.

Choose an isomorphism of $k$-algebras $O(\Gr) \stackrel{\sim}{\to} k((\ul \fg^{(r)})_{(1)}\times \bG_{(r)})$.
Then sending a $\pi$-point pair $(\alpha_K,\beta_K)$ of $\tilde D(\bG_{(r)})$ to 
$$\alpha_K + \beta_K: K[t]/t^p \ \to \ O(\bG_{(r),K}) \ \stackrel{\sim}{\to} \  K( (\ul \fg^{(r)})_{(1)}\times \bG_{(r)})$$
determines a bijection 
\begin{equation}
\label{eqn:bij}
\Pi_O: \Pi (O(\bG_{(r)})) \stackrel{\sim}{\to} \  \Pi ((\ul \fg^{(r)})_{(1)} \times \bG_{(r)} ),
\end{equation}
where  we abuse notation by interpreting the right hand side as the set of (Zariski) points of the 
indicated $\Pi$-point scheme.  

Moreover,
for every $O(\bG_{(r)})$-module $M$, $\Pi_O$ restricts to a bijection
\begin{equation}
\label{eqn:bijM}
\Pi_{O,M}: \Pi (O(\bG_{(r)}))_{M} \  \simeq \ \Pi((\ul \fg^{(r)})_{(1)} \times \bG_{(r)})_M.
\end{equation}
\end{prop}

\begin{proof}
Since $(\ul \fg^{(r)})_{(1)}$ is an abelian unipotent finite group scheme,  since $\beta_K$ is a $\pi$-point
of $\bG_{(r)}$ and thus factors through an abelian unipotent subgroup scheme, and since $\alpha_K + \beta_K: 
K[t]/t^p \to KG$ is flat by Lemma \ref{lem:flat}, we conclude 
that $\alpha_K + \beta_K$ is a $\pi$-point of the finite group scheme $(\ul \fg^{(r)})_{(1)} \times \bG_{(r)}$.
Moreover, the discussion preceding shows that every point of $\Pi ((\ul \fg^{(r)})_{(1)} \times \bG_{(r)})$
(i.e., every equivalence class of $\pi$-points of $\ul \fg^{(r)})_{(1)} \times \bG_{(r)}$) is represented by such
a flat map $\alpha_K + \beta_K$ associated to a $\pi$-point pair of $\tilde D(\Gr)$.  

Thus, to prove that
$\Pi_O$ of (\ref{eqn:bij}) is well defined and injective (as well as surjective), it suffices to show that
two $\pi$-point pairs $(\alpha_K,\beta_K), \ (\alpha_K^\prime,\beta_K^\prime)$ of $\tilde D(\Gr)$
are $O(\Gr)$-equivalent if and only if $\alpha_K+\beta_K, \ \alpha_K^\prime + \beta_K^\prime$ are 
equivalent $\pi$-points of $(\ul \fg^{(r)})_{(1)} \times \bG_{(r)}$.  This follows immediately from
Definition \ref{defn:equiv}, Proposition \ref{prop:equal-prod}, and the bijection 
$\Psi_{\ul \fg^{(r)})_{(1)} \times \bG_{(r)}}$ of (\ref{eqn:PiG}).

The fact that (\ref{eqn:bij}) restricts to a bijection $\Pi (O(\bG_{(r)}))_{M} \ 
\stackrel{\sim}{\to} \  \ \Pi ((\ul \fg^{(r)})_{(1)} \times \bG_{(r)})_{M}$ follows immediately from the definitions
(i.e., from Definition \ref{defn:equiv} plus Proposition \ref{prop:Pi-equiv} and \cite[Defn 3.2]{F-P1}).
\end{proof}

\vskip .2in


\section{The tensor product property for $\tilde D(\Gr)$-modules}
\label{sec:tensor}

Propositions \ref{prop:Pi-equiv} and  \ref{prop:OtoG} identify $M \mapsto \Pi(\tilde D(\Gr))_M$ with $M \mapsto \Pi((\ul \fg^{(r)})_{(1)} \times \bG_{(r)})_M$.
This immediately tells us that $M \mapsto \Pi(\tilde D(\Gr)_M$ satisfies most of the properties required for a good
theory of support varieties of $\tilde D(\Gr)$-modules.  However, this identification does not imply that
$M \mapsto \Pi(\tilde D(\Gr)_M$ satisfies the ``tensor product property" (namely, that 
$ \Pi(\tilde D(\Gr)_{M\otimes M^\prime} \ = \ \Pi(\tilde D(\Gr)_M \cap \Pi(\tilde D(\Gr)_{M^\prime}$) 
because $O(\bG_{(r),K}) \ \stackrel{\sim}{\to} \  K( (\ul \fg^{(r)})_{(1)}\times \bG_{(r)})$
is not a map of Hopf algebras.   In this section, we prove this tensor product property.  
 
The following proposition of \cite{F-P1} is the key to our proof of 
Theorem \ref{thm:2coproducts} establishing the tensor product property for $M \mapsto \Pi(\tilde D(\Gr))_M$.

\begin{prop} \cite[Prop 2.2]{F-P1}
\label{prop:FP1}
Let $V$ be a $k$-vector space and $\alpha, \beta, \gamma$ be pairwise commuting $V$-endomorphisms.
Assume further $\alpha, \beta$ are $p$-nilpotent and $\gamma$ is $p^r$-nilpotent in $End_k(V)$
for some $ \geq 1$.  Then $V$ is 
projective as a $k[u]/u^p$-module where the action of $u$ is given by $\alpha$ if and only if it is projective as a
$k[v]/v^p$module where the action of $v$ is given by $\alpha + \beta\gamma$.
\end{prop}

We apply Proposition \ref{prop:FP1} to prove the following comparison of pull-backs along a given $\pi$-point
of tensor products of $A$-modules with two different $A$-module structures determined by different 
coproducts on $A$.

\begin{prop}
\label{prop:pullback}
Let $A$ be a finite dimensional, local, commutative $k$-algebra with maximal ideal $\fm$ 
satisfying the property that
$a^p = 0$ for all $a \in \fm$.  Assume that $A$ is equipped with  two $A$-linear coproducts
$\Delta, \ \Delta^\prime: A \to A\otimes A$ with the property that
\begin{equation}
\label{eqn:Delta-form}
(\Delta(a) - a\otimes 1 - 1 \otimes a), \ (\Delta^\prime(a) - a\otimes 1 - 1 \otimes a) \ \in \ \fm \otimes \fm,
\quad \forall a \in \fm.
\end{equation}
Let $M, \ M^\prime$ be $A$-modules and let $M \boxtimes M^\prime$ be the $A\otimes A$ module 
given as the external tensor
product.  Then for any  flat map $\alpha: k[t]/t^p \to A$, the pull-back of $M \boxtimes M^\prime$ via
$\Delta \circ \alpha$ is a projective $k[t]/t^p$-module if and only if the pull-back  of $M \boxtimes M^\prime$ via
$\Delta^\prime \circ \alpha$ is a projective $k[t]/t^p$-module. 

More generally, let $H$ be  finite dimensional, local, commutative Hopf algebra with maximal ideal $\fm_H$ 
such that $c^p = 0$ for all $c \in \fm_H$ and equip $A\otimes H$ with the coproducts 
$\Delta \otimes \Delta_H, \ \Delta^\prime \otimes \Delta_H$.  
Consider $(A\otimes H)$-modules $M, \ M^\prime$, and let $M \boxtimes M^\prime$ denote the
$(A\otimes H)\otimes (A\otimes H)$ module given as the external tensor product.   Then for any
flat map $\alpha: k[t]/t^\p \to A\otimes H$, the pull-back of $M \boxtimes M^\prime$ via
$(\Delta \otimes \Delta_H )\circ \alpha$ is a projective $k[t]/t^p$ module  if and only if the pull-back  of
 $M \boxtimes M^\prime$ via $(\Delta^\prime \otimes \Delta_H) \circ \alpha$ is a projective $k[t]/t^p$-module. 
\end{prop}

\begin{proof}
We apply Proposition \ref{prop:FP1} to the $k$-vector space $V = 
M\otimes M^\prime$ and consider triples of $p$-nilpotent elements in the image of $A\otimes A \to
End_k(V)$ given by viewing $V$ as the external tensor product $M\boxtimes M^\prime$ 
(and thus an $A\otimes A$ module for the commutative $k$-algebra $A\otimes A$).

Consider an element $a \in \fm \subset A$ and write $\Delta(a) - \Delta^\prime(a) \  = \
\sum_{i=1}^t \ b_i \otimes b_i^\prime \in \fm \otimes \fm.$  Observe that 
$\Delta(a) - \Delta^\prime(a) \in A\otimes A$ acts on $V = M\boxtimes M^\prime$ as the sum of the actions 
commuting elements $(b_i\otimes 1)\cdot (1\otimes b_i^\prime), 1 \leq i \leq ,t$ each of which is $p$-nilpotent..  
We apply Proposition \ref{prop:FP1} with $V = M\boxtimes M^\prime$ successively to the triples 
$$
\alpha = ((\Delta^\prime(a) + \sum_{s=1}^i b_s \otimes b_s^\prime), \quad \beta = (b_{i+1} \otimes 1),
\quad \gamma =  (1\otimes b_{i+1}^\prime)$$ 
with $i+1 \leq t$ to conclude
that $V$ is projective for $k[u]/u^p$ with $u$ acting as  $\Delta^\prime(a) + \sum_{s=1}^i b_s \otimes b_s^\prime$
if and only if $V$ is projective for $k[v]/t^p$ with $v$ acting as 
$\Delta^\prime(a) + \sum_{s=1}^{i+1} b_s \otimes b_s^\prime$.  This establishes the first assertion.

We extend this argument to $A \otimes H$, a commutative algebra with coproducts $\Delta \otimes \Delta_H, \ 
\Delta^\prime \otimes \Delta_H$.  We denote by $\fn \subset A \otimes H$ the maximal ideal of $A \otimes H$
generated by $\fm \otimes 1, \ 1 \otimes \fm_H$,
and observe that every element in $\fn$ has $p$-th power 0.  Moreover, 
$(\Delta \otimes \Delta_H) - (\Delta^\prime \otimes \Delta_H) \ = \ (\Delta - \Delta^\prime) \otimes \Delta_H$ 
applied to $a\otimes 1 \in \fn$ lies in  $\fn \otimes \fn$ whereas $(\Delta - \Delta^\prime) \otimes \Delta_H$ 
applied to $1\otimes h \in 1 \otimes \fm_H$ is 0.  Consequently, 
$$((\Delta \otimes \Delta_H) - (\Delta^\prime \otimes \Delta_H))(b) \ \in \ \fn\otimes \fn, \quad \forall b \in \fn.$$
Thus we may repeat the preceding argument with $(A,\fm)$ replaced by $(A\otimes H,\fn)$ to complete 
the proof.
\end{proof}

\begin{defn}
\label{defn:other-tensor}
Consider $O(\Gr)$-modules $M, \ M^\prime$.

We denote by $M \otimes M^\prime$ the $O(\Gr)$-module determined by the coproduct of $O(\Gr)$
defined as  by restriction of the coproduct of $\tilde D(\Gr)$ (see Proposition \ref{prop:OG}).

We denote by $M\otimes_G M^\prime$ the $k$-vector space  $M \otimes M^\prime$ equipped
with the $O(\Gr) \ \simeq \ (k((\ul \fg^{(r)})_{(1)}) \times \bG_{(r)})$-module structure determined 
by the coproduct for the group algebra $k((\ul \fg^{(r)})_{(1)}) \times \bG_{(r)})$ of the finite
group scheme $(\ul \fg^{(r)})_{(1)}) \times \bG_{(r)}$. 
\end{defn}

In the following theorem, we prove the ``tensor product property" for $\tilde D(\Gr)$-modules; namely,
$M \ \mapsto \ \Pi(\tilde D(\bG_{(r)}))_M$ sends tensor product of modules to intersection of the supports
of its tensor factors.

\begin{thm}
\label{thm:2coproducts}
Let $(\alpha_K,\beta_K)$ be a  $\pi$-point pair of $\tilde D(\Gr)$ as in Definition \ref{defn:pairs} and let 
$M, \ M^\prime$ be $O(\Gr)$-modules.  Then $(\alpha_K+\beta_K)^*(M_K\otimes M^\prime_K)$ is a free
$K[t]/t^p$-module if and only if $(\alpha_K+\beta_K)^*(M_K\otimes_G M^\prime_K)$
is a free $K[t]/t^p$-module.  Thus,
\begin{equation}
\label{eqn:tO}
\Pi(\tilde D(\bG_{(r)}))_{M\otimes M^\prime}  \ = \ \Pi(\tilde D(\bG_{(r)}))_M \cap \Pi(\tilde D(\bG_{(r)}))_{M^\prime}.
\end{equation}
\end{thm}
 
\begin{proof}
We recall that any $\beta_K: K[t]/t^p \to K\Gr$ is equivalent to a map which factors through 
a map of group algebras induced by a map of finite group schemes
$ \bG_{a(r),L} \to \bG_{(r),L}$ \cite[Prop 4.2]{F-P1}.  By replacing $K$ if 
necessary by a field extension $K^\prime/K$, we may assume $K = L$.   Thus when comparing
restrictions along a given $\pi$-point pair, we may restrict to 
$K[\bG_{(r+1)}/\Gr] \otimes K\bG_{a(r)} \hookrightarrow O(\Gr)\otimes K$
and replace $M, \ M^\prime$ by their restrictions to $K[\bG_{(r+1)}/\Gr] \otimes K\bG_{a(r)}$.

We apply Proposition \ref{prop:pullback} with $A$ equal to  $K[\bG_{(r+1)}/\Gr]$ and $H = K\bG_{a(r)}$,
and equip $A\otimes K\bG_{a(r)}$ with the two coproduct structures obtained by restriction of the two coproduct structures 
on $O(\Gr)$ discussed in Definition \ref{defn:other-tensor}.   As shown in \cite[I.2.4]{J}, $\Delta, \  \Delta^\prime$
satisfy condition (\ref{eqn:Delta-form}) of Proposition \ref{prop:pullback}.   In our application of Proposition \ref{prop:pullback},
we take $\alpha:K[t]/t^p \ \to \  K[\bG_{(r+1)}/\Gr] \otimes K\bG_{a(r)}$
to be the composition of the coproduct $K[t]/t^p \to K[t]/t^p\otimes K[t]/t^p, \ t \mapsto t\otimes 1 + 1 \otimes t$ 
and the external tensor product $\alpha_K \otimes \beta_K: K[t]/t^p \to K[t]/t^p\otimes K[t]/t^p \to K[\bG_{(r+1)}/\Gr]\otimes K\bG_{a(r)}$.

By applying Proposition \ref{prop:pullback} with these values for $A, H, \alpha$, we conclude that it suffices
to prove the equality (\ref{eqn:tO}) after replacing the tensor product for \\
$K[\bG_{(r+1)}/\Gr] \otimes K\bG_{a(r)}$-modules
(given by the restriction of the coproduct of $O(\Gr)$) by the tensor product given by the restriction 
of the tensor product
of $(\fg^{(r)})_{(1)}\times \bG_{(r)}$-modules.   To conclude the proof of  equality (\ref{eqn:tO}), 
we apply the tensor product property for modules for the 
(infinitesimal) group scheme $G \ = \ (\fg^{(r)})_{(1)}\times \bG_{(r)}$ \cite[Prop 3.2]{F-P2}; in other words, for $M\otimes_G M^\prime$
in the notation of Definition \ref{defn:other-tensor}.
\end{proof}

Theorem \ref{thm:2coproducts} enables topological structures on the sets $\Pi(\tilde D(\Gr))$ and $\Pi(O(\Gr))$.

\begin{cor}
\label{cor:topologies}
Defining a subset of $\Pi(\tilde D(\Gr))$ to be closed if and only if it is 
of the form $\Pi(\tilde D(\Gr))_M \subset \Pi(\tilde D(\Gr))$ for some finite dimensional
$\tilde D(\Gr)$-module $M$ defines a topology on $\Pi(\tilde D(\Gr))$.
Similarly, defining a subset of $\Pi(O(\Gr))$ to be closed if and only if it is 
of the form $\Pi(O(\Gr))_M \subset \Pi(O(\Gr))$ for some finite dimensional
$O(\Gr)$-module $M$ defines a topology on $\Pi(O(\Gr))$.

With these topologies, the bijection $\Pi_{\tilde D}: \Pi (\tilde D(\bG_{(r)})) \ \twoheadrightarrow \ \Pi (O(\bG_{(r)}))$
of Proposition \ref{prop:Pi-equiv} is a homeomorphism.
\end{cor}

\begin{proof}
The asserted formulation of a topology on $\Pi(\tilde D(\Gr))$ is justified by Proposition \ref{prop:basic}(4) and
Theorem \ref{thm:2coproducts}.  This applies as well to justifiy the topology on $\Pi(O(\Gr))$, since 
$i_O: O(\Gr) \to \tilde D(\Gr)$ is a map of Hopf algebras.
\end{proof}

\vskip .2in


\section{$\Psi_{\tilde D}: \Pi(\tilde D(\Gr)) \to \bP(H^\bu(\tilde D(\Gr),k))$}
\label{sec:cohom}

Cohomological support varieties can be readily defined for any finite dimension Hopf algebra
$H$ whose cohomology is finitely generated (i.e., such that $H^*(H,k)$ is a finitely generated 
$k$-algebra and that $H^*(H,M)$ is a finite $H^*(H,k)$-module for any finite dimensional $H$-module). 
 In this section, we exhibit in Proposition \ref{prop:cont} a continuous map $\Psi_{\tilde D}$ 
from the $\pi$-point pair space $\Pi(\tilde D(\Gr))$ introduced in Section \ref{sec:pairs}
to the (projectivized) cohomology support variety $\bP(H^\bu(\tilde D(\Gr),k))$. To do this,
we utilize the subalgebra $O(\Gr) \subset \tilde D(\Gr)$ and the related finite group scheme 
$\Gr \times (\ul\fg^{(r)})_{(1)}$, relying upon the isomorphism $\Pi(\Gr \times (\ul\fg^{(r)})_{(1)}) \
\stackrel{\sim}{\to} \ \bP(H^\bu(\Gr \times (\ul\fg^{(r)})_{(1)},k)$.

Associated to a 
finite dimensional $\tilde D(\Gr)$-module $M$, we consider the annihilator ideal 
$I_M \subset H^\bu(\tilde D(\Gr),k)$ of the $H^\bu(\tilde D(\Gr),k)$-module $Ext_{\tilde D(\Gr)}^*(M,M)$ 
which equals the annihilator of $1 \in Ext_{\tilde D(\Gr)}^*(M,M)$ which in turn equals the kernel 
of the graded map of $k$-algebras $H^\bu(\tilde D(\Gr),k) \to Ext_{\tilde D(\Gr)}^*(M,M)$.
We shall consider the projective scheme $\bP(H^\bu(\tilde D(\Gr),k))$  whose (Zariski) points 
are  (non-trivial) homogenous prime ideals of $H^\bu(\tilde D(\Gr),k)$ and the reduced projective scheme 
$\bP(H^\bu(\tilde D(\Gr),k))_M$ whose points are  (non-trivial) homogenous prime ideals of 
$H^\bu(\tilde D(\Gr),k)$ containing $I_M$.

Consider the quotient map $\bG_{(r+2)} \ \twoheadrightarrow \ (\bG^{(r+1)})_{(1)}$, the image 
of $F^{r+1}: \bG_{(r+2)} \to (\bG^{(r+1)})_{(r+2)}$ whose kernel is $\bG_{(r+1)}$.  This quotient map
determines the flat map of $k$-algebras
$$k[(\bG^{(r+1)})_{(1)}] \ \simeq \ k[\bG_{(r+2})/\bG_{(r+1)}] \ \to \ k[\bG_{(r+2)}]$$
with the property that $k[\bG_{(r+1)}] \ \simeq \ k[\bG_{(r+2)}] \otimes_{k[\bG_{(r+2})/\bG_{(r+1)}]}k$.

The following proposition is a summary of a central technique of \cite{F-N}.

\begin{prop}
\label{prop:theta}
The deformation $k[\bG_{(r+2})/\bG_{(r+1)}] \ \to \ k[\bG_{(r+2)}]$ of $k[\bG_{(r+1)}]$ (parametrized by
$k[\bG_{(r+2})/\bG_{(r+1)}] \simeq k[(\bG^{(r+1)})_{(1)}]$) determines a $\bG$-equivariant map 
\begin{equation}
\label{eqn:Gerst1}
\fg^{(r+1)} \ \to \ H^2(k[\bG_{(r+1)}],k).
\end{equation}   

The deformation $k[\bG_{(r+2})/\bG_{(r+1)}] \# k\Gr \ \to \ k[\bG_{(r+2)}] \# k\Gr$ (also parametrized 
by $k[\bG_{(r+2})/\bG_{(r+1)}]$) provides  a $\Gr$-equivariant lifting of (\ref{eqn:Gerst1})
\begin{equation}
\label{eqn:Gerst2}
\sigma_{\tilde D}: \fg^{(r+1)} \ \to \ H^2(\tilde D(\Gr),k).
\end{equation}.

Consequently, we obtain a map of $k$-algebras
\begin{equation}
\label{eqn:t-theta}
\tilde \theta_r: H^*(\bG_{(r)},k) \otimes S^\bu(\fg^{(r+1)}[2]) \ \to \ H^*(\tilde D(\bG_{(r)}),k)
\end{equation}
defined as the product of the inflation map on $H^*(\bG_{(r)},k)$ and the algebra extension of $\sigma_{\tilde D}$
on $S^\bu(\fg^{(r+1)}[2])$.
\end{prop}

\begin{proof}
Each element of the tangent space at the identity of $(\bG^{(r+1)})_{(1)}$ (naturally identified with $\fg^{(r+1)}$)
corresponds to a map $k[(\bG^{(r+1)})_{(1)}] \to k[\epsilon]/\epsilon^2$ of $k$-algebras.  
Extending the deformation $k[\bG_{(r+2})/\bG_{(r+1)}] \ \to \ k[\bG_{(r+2)}]$ along such a map determines
a deformation of $k[\bG_{(r+1)}]$ parametrized by $k[\epsilon]/\epsilon^2$.  Thus, the map
(\ref{eqn:Gerst1}) is a consequence of 
the naturality of Gerstenhaber's bijection  \cite{Ger} between deformations of $k[\bG_{(r+1)}]$ parametrized 
by $k[\epsilon]/\epsilon^2$ and the second Hochschild cohomology group $HH^2(k[\bG_{(r+1)}])$ of $k[\bG_{(r+1)}]$
together with the natural map $HH^2(k[\bG_{(r+1)}]) \to H^2(k[\bG_{(r+1)}],k))$ (see \cite[Prop 3.4]{F-N}).
The map $\sigma_{\tilde D}$ of (\ref{eqn:Gerst2}) follows exactly the same way.  
Moreover, the naturality of Gerstenhaber's bijection implies that
$\sigma_{\tilde D}$ lifts  (\ref{eqn:Gerst1}).

The construction of $\tilde \theta_r: H^*(\bG_{(r)},k) \otimes S^\bu(\fg^{(r+1)}[2]) \ \to \ H^*(\tilde D(\bG_{(r)}),k)$ is 
implicit in the last statement of this proposition.
\end{proof}

Finite generation of the cohomology algebra $H^*(D(\bG_{(r)}),k)$ and the finiteness of $H^*(D(\bG_{(r)}),M)$ as 
an $H^*(D(\bG_{(r)}),k)$-module  (for a finite dimensional $D(\Gr)$-module $M$) was established in \cite[Thm 5.3]{F-N}.
Using Proposition \ref{prop:theta}, we extend this finite generation property to $\tilde D(\Gr)$ and $O(\Gr)$ as stated below;
the proof  is a straight-forward adaptation of that of \cite[Thm 5.3]{F-N} which in turn is an adaptation of \cite[Thm 1.1]{F-S}.

\begin{thm}
\label{thm:fingen}
Consider the map \ $\tilde \theta_r: H^*(\bG_{(r)},k) \otimes S^\bu(\fg^{(r+1)}[2]) \ \to \ H^*(\tilde D(\bG_{(r)}),k)$
of (\ref{eqn:t-theta}).
\begin{enumerate}
\item
The composition 
\begin{equation}
i_O^* \circ \tilde \theta_r: H^*(\bG_{(r)},k) \otimes S^\bu(\fg^{(r+1)}[2]) \ \to \ H^*(\tilde D(\bG_{(r)}),k)
\end{equation} 
$$\to \ H^*(O(\bG_{(r)}),k) \ \simeq \ H^\bu(\bG_{(r)},k) \otimes H^*(k[(\bG^{(r)})_{(1)}],k)$$
is a split injection with $H^*(\bG_{(r)},k) \otimes S^\bu(\fg^{(r+1)}[2])$-module complement
$$H^*(\bG_{(r)},k) \otimes \Lambda^{> 0}(\fg[1]) \otimes  \ S^\bu(\fg^{(r+1)}[2])$$ 
which consists entirely of nilpotent elements 
\item
For any finite dimensional $\tilde D(\bG_{(r)})$-module $M$,
$H^*(\tilde D(\bG_{(r)}),M)$ is a finitely generated 
$H^\bu(\bG_{(r)},k) \otimes S^\bu(\fg^{(r+1)}[2])$-module (with action given by the
restriction along $\tilde \theta_r$ of the natural action of .$H^*(\tilde D(\bG_{(r)}),k)$ on $H^*(\tilde D(\bG_{(r)}),M)$).
\item
For any finite dimensional $O(\Gr)$-module $M$, $H^*(O(\bG_{(r)}),M)$ is also a 
finitely generated $H^\bu(\bG_{(r)},k) \otimes S^\bu(\fg^{(r+1)}[2])$-module.
\end{enumerate}
\end{thm}

\begin{proof}
We utilize the spectral sequence 
\begin{equation}
\label{eqn:Grothss}
E_2^{s,t} = H^s(\bG_{(r)}, H^t(k[G_{(r+1)}],k)) \ \Rightarrow \ H^{s+t}(\tilde D(\bG_{(r)}),k),
\end{equation}
which can be derived as a special case of Grothendieck's spectral sequence for the derived 
functors of a composition of left exact functors exactly as in \cite[Prop 5.2]{F-N}. 

The composition $i_O^* \circ \tilde \theta_r$ restricts to the identity on the tensor factor $H^*(\bG_{(r)},k)$
by the definition of $\tilde \theta_r$ and the fact that the composition $k\Gr  \to O(\Gr) \to \tilde D(\Gr)$ is the evident
inclusion. The deformation construction leading to the definition of $\tilde \theta_r$ is designed to insure that
the composition $S^\bu(\fg^{(r+1)}[2]) \to H^*(\tilde D(\bG_{(r)}),k) \to H^*(k[\bG_{(r+1)}/\Gr],k)$ sends elements
of $\fg^{(r+1)}[2]$ to elements of $H^2(k[\bG_{(r+1)}/\Gr],k)$ which generate $H^*(k[\bG_{(r+1)}/\Gr],k)$
modulo nilpotent elements.   The  identification of the restriction of the map
$i_O^* \circ \tilde \theta_r$ to $1\otimes S^\bu(\fg^{(r+1)}[2])$ and the computation of $H^*(k[(\bG^{(r)})_{(1)}],k)$  
are given in \cite[\S 4]{F-N}, thereby 
establishing the asserted split surjection.

The fact that $H^*(\tilde D(\bG_{(r)}),M)$ is a finitely generated $H^\bu(\bG_{(r)},k) \otimes S^\bu(\fg^{(r+1)}[2])$-module 
for any finite dimensional $\tilde D(\bG_{(r)})$-module $M$ arises from the natural pairing 
of  (\ref{eqn:Grothss}) with the corresponding spectral sequence
\begin{equation}
\label{eqn:GrothM}
E_2^{s,t}(M) = H^s(\bG_{(r)}, H^t(k[G_{(r+1)}],M)) \ \Rightarrow \ H^{s+t}(\tilde D(\bG_{(r)}),M)
\end{equation}
with $M$ coefficients (as in the proof of \cite[Prop 5.2]{F-N}).  The same argument with 
$\tilde D(\bG_{(r)})$ replaced by $O(\Gr)$ implies that  $H^*(O(\bG_{(r)}),M)$ is a finitely generated 
$H^\bu(\bG_{(r)},k) \otimes S^\bu(\fg^{(r+1)}[2])$-module.
\end{proof}

\begin{cor}
\label{cor:finite}
The  map of commutative $k$-algebras $\tilde \theta_r$ induces a finite, surjective map 
$$\tilde \Theta_r: \Spec H^\bu(\tilde D(\bG_{(r)}),k) \to \Spec H^\bu(\Gr,k) \times \ul{\fg^*}^{(r+1)}.$$

Consequently,
$$ dim (\Spec H^\bu(\tilde D(\bG_{(r)}),k)) \ = \ dim(\Spec H^\bu(\Gr,k) + dim(\fg).$$
\end{cor}

\begin{proof}
By Theorem \ref{thm:fingen}, the composition $i^* \circ \tilde \theta_r$ is injective, so that $\tilde \theta_r$
is injective; this theorem also tells us that $\tilde \theta_r$ is a finite map.  Consequently,
$\tilde \theta_r$ is an integral extension which implies that $\Theta_r$ is finite and surjective (see \cite[Thm 9.3]{Mat}).

The determination of  $dim(\Spec H^\bu(\tilde D( \bG_{(r)}),k))$ follows from the observation that 
$dim( \Spec H^\bu(\Gr,k) \times \ul{\fg^*}^{(r+1)} ) \ =  \ dim(\Spec H^\bu(\Gr,k)) + dim(\fg).$
\end{proof}

We establish the following notation.

\begin{note}
\label{note:Proj}
We denote by $\bP (H^\bu(\tilde D(\bG_{(r)}),k))$ (respectively, $\bP (H^\bu(O(\bG_{(r)}),k))$; respectively, 
$\bP (H^\bu(\Gr,k))$)
the projective scheme associated to the commutative graded $k$-algebra $H^\bu(\tilde D(\bG_{(r)}),k)$
(resp., $H^\bu(O(\bG_{(r)}),k)$;  resp., $H^\bu(\bG_{(r)},k))$).  Thus, 
the points of $\bP (H^\bu(\tilde D(\bG_{(r)}),k))$ (resp., 
$\bP (H^\bu(O(\bG_{(r)}),k)$; resp., $\bP (H^\bu(\bG_{(r)},k))$) are the non-trivial homogeneous prime ideals of 
$H^\bu(\tilde D(\bG_{(r)}),k)$ (resp., $H^\bu(O(\bG_{(r)}),k)$; resp., $H^\bu(\bG_{(r)},k))$).  
For any finite dimensional $\tilde D(\bG_{(r)})$-module $M$, we denote by
 $$ \bP (H^\bu(\tilde D(\bG_{(r)}),k))_M \ \subset \ \bP (H^\bu(\tilde D(\bG_{(r)}),k)),$$
 $$ \bP (H^\bu(O(\bG_{(r)})),k))_M \ \subset \ \bP (H^\bu(O(\bG_{(r)})),k)),$$
 $$ \bP (H^\bu(\bG_{(r)},k))_M \ \subset \ \bP (H^\bu(\bG_{(r)},k))$$
 the reduced closed subscheme given by the radicals of the ideals 
 $$ker\{ H^\bu(\tilde D(\bG_{(r)}),k) \to Ext^*_{\tilde D(\bG_{(r)})}(M,M)\},$$
 $$ker\{ H^\bu(O(\bG_{(r)}),k) \to Ext^*_{O(\bG_{(r)}))}(M,M)\},$$
 $$ker\{ H^\bu(\bG_{(r)},k) \to Ext^*_{\bG_{(r)}}(M,M)\}$$
respectively.
\end{note}

\vskip .1in

The equivalence relation on $\pi$-point pairs $(\alpha_K,\beta_K)$ of $\tilde D(\Gr)$ (and thus on $O(\Gr)$)
 has been designed to enable the validity of the following proposition.

\begin{prop}
\label{prop:PsiO}
Define $\Psi_O: \Pi(O(\Gr)) \to \bP(H^*(O(\Gr),k))$ to be the map of sets sending the equivalence class of 
$(\alpha_K,\beta_K) \in \Pi(O(\Gr))$ to 
$$ker\{ (\alpha_K+\beta_K)^*: H^\bu(O(\Gr)_K,K) \to H^\bu(K[t]/t^p,K) \} \ \cap \ H^*(O(\Gr),k)$$
in $\bP(H^*(O(\Gr),k))$.
A choice of isomorphism of $k$ algebras $j: O(\Gr) \ \stackrel{\sim}{\to} \ k(\ul \fg^{(r)})_{(1)} \times \bG_{(r)})$
as in Proposition \ref{prop:ulg} determines the vertical maps in the following commutative 
square of homeomorphisms
\begin{equation}
\label{comm:PsiO}
\begin{xy}*!C\xybox{%
\xymatrix{
\Pi(O(\Gr))\ar[r]^-{\Psi_O} \ar[d]_{\Pi_O} & \bP(H^\bu(O(\Gr),k)) \ar[d]^{\bP_O} \\
\Pi(\ul \fg^{(r)})_{(1)} \times \bG_{(r)}) \ar[r]^-{\Psi_G} & \bP(H^\bu(\ul \fg^{(r)})_{(1)} \times \bG_{(r)},k)) }
}\end{xy}
\end{equation}
whose lower horizontal map is the natural homeomorphism of C for $G$  the finite group scheme
$(\ul \fg^{(r)})_{(1)} \times \bG_{(r)}$.

Moreover, for every finite dimensional $O(\Gr)$-module $M$, $\Psi_O$ restricts to a bijection
$\Psi_{O,M}: \Pi(O(\Gr))_M \ \stackrel{\sim}{\to} \ \bP(H^\bu(O(\Gr),k)) _M$.
\end{prop}

\begin{proof}   The definition of $\Psi_O$  is compatible with the definition of 
$\Psi_G: \Pi(\ul \fg^{(r)})_{(1)} \times \bG_{(r)}) \ \to \  \bP(H^\bu(\ul \fg^{(r)})_{(1)} \times \bG_{(r)},k))$
(see $(\ul \fg^{(r)})_{(1)}\times \bG_{(r)}$ and \cite[Prop 3.5]{F-P1}), so that (\ref{comm:PsiO}) commutes.  
The fact that $\Psi_O$ is well defined 
follows from the fact that $\Psi_G$ is well defined and that the maps $\Pi_O, \ \bP_O$ are bijections.  
The topology on
$\Pi(\tilde D(\Gr)) \ = \ \Pi(O(\Gr))$ introduced in Theorem \ref{thm:2coproducts} is equal to that induced 
by the bijection $\Pi_O$,
so that this bijection is a homeomorphism.  The bijection $\bP_O$ is a homeomorphism since the isomorphism $j$
induces an isomorphism between the category of finite dimensional $O(\Gr)$-modules and the category of finite 
dimensional $k(\ul \fg^{(r)})_{(1)} \times \bG_{(r)})$-modules.

The fact that $\Psi_{O,M}$ is a bijection for any finite dimensional $O(\Gr)$-module follows 
from the fact that $\Psi_O$ restricts to a bijection
$\Psi_{O,M}: \Pi(\ul \fg^{(r)})_{(1)} \times \bG_{(r)})  \ \stackrel{\sim}{\to} \ 
\bP(H^\bu(\ul \fg^{(r)})_{(1)} \times \bG_{(r)},k))_M$.
\end{proof}

The following proposition complements Proposition \ref{prop:PsiO} while introducing 
$\Psi_{\tilde D}$.

\begin{prop}
\label{prop:cont}
The composition $\Psi_{\tilde D} \equiv \bP_{\tilde D} \circ \Psi_O \circ \Pi_{\tilde D}^{-1}$:  
$$\Pi(\tilde D(\Gr)) \stackrel{\sim}{\to} \Pi(O(\Gr)) \ \stackrel{\sim}{\to} \ \bP(H^\bu(O(\Gr),k))
 \ \to \ \bP(H^\bu(\tilde D(\Gr),k))$$
is continuous, fitting in the commutative square
\begin{equation}
\label{comm:PiO}
\begin{xy}*!C\xybox{%
\xymatrix{
\Pi(\tilde D(\Gr))\ar[r]^-{\Psi_{\tilde D}} & \bP(H^\bu(\tilde D(\Gr),k)) \\
\Pi(O(\Gr))  \ar[u]_{\Pi_{\tilde D}}^-{\simeq}   \ar[r]_-{\Psi_O}^-{\sim} & \bP(H^\bu(O(\Gr),k)) \ar[u]^{\bP_{\tilde D}}}
}\end{xy}
\end{equation}
with the properties that $\Pi_{\tilde D}$ and $\Psi_O$ are homeomorphisms, and that the map
$\bP_{\tilde D}$ (induced by $i_O: O(\Gr) \to \tilde D(\Gr)$) is a closed immersion.

For any finite dimensional $\tilde D(\Gr)$-module $M$, $\Psi_{\tilde D}$ restricts to
\begin{equation}
\label{eqn:PsiDM}
\Psi_{\tilde D,M}: \Pi(\tilde D)_M \ \hookrightarrow \ \bP(H^\bu(\tilde D(\Gr),k))_M.
\end{equation}
\end{prop}

\begin{proof} 
By Proposition \ref{prop:Pi-equiv}, $\Pi_{\tilde D}$ is a homeomorphism, whereas
 $\Psi_O$ is a homeomorphism by Proposition \ref{prop:PsiO}. 
Theorem \ref{thm:fingen}(1) implies that $i_O^*$ is a finite map and is surjective modulo nilpotents
which implies that $\bP_O$ is a closed embedding.  Thus, $\Psi_{\tilde D}$ is also a closed
embedding and thus is continuous.

To prove that $\Psi_{\tilde D}(\Pi(\tilde D(\Gr))_M)$ is contained in $\bP(H^\bu(\tilde D(\Gr),k))_M$, 
we consider some
$(\alpha_K,\beta_K)$ representing a point of $\Pi (\tilde D(\bG_{(r)}))_M$; in other words, 
assume that $(\alpha_K+\beta_K)^*(M_K)$ is not projective.  Consider the following commutative square
\begin{equation}
\label{comm:Psi_M}
\begin{xy}*!C\xybox{%
\xymatrix{
H^\bu(\tilde D(\bG_{(r)}),k) \ar[r]^-{(\alpha_K+\beta_K)^*} \ar[d] & H^\bu(K[t]/t^p,K) \ar[d] \\
Ext^*_{\tilde D(\bG_{(r)}))}(M,M) \ar[r] & Ext^*_{K[t]/t^p}((\alpha_K+\beta_K)^*(M_K),(\alpha_K+\beta_K)^*(M_K))}
}\end{xy}
\end{equation}
whose horizontal maps are induced by restriction along $\alpha_K + \beta_K$.
Since  $(\alpha_K+\beta_K)^*(M_K)$ is not projective, the right vertical map is injective.
A simple diagram chase  then implies that $ker \{ (\alpha_K+\beta_K)^*: H^\bu(\tilde D(\bG_{(r)})),k) \to H^\bu(K[t]/t^p,K)\}$
contains $ker\{ H^\bu(\tilde D(\bG_{(r)})),k) \to Ext^*_{K[t]/t^p}((\alpha_K,\beta_K)^*(M_K),(\alpha_K,\beta_K)^*(M_K))\}$, 
so that $\Psi_{\tilde D}$ applied to the class of  $(\alpha_K,\beta_K)$ lies in $\bP (H^\bu(\tilde D(\bG_{(r)}),k))_M$.
\end{proof}

\vskip .2in


\section{$\Pi(\tilde D(\Gr))$ for $\bG$ admitting a quasilogarithm}
\label{sec:quasi-log}

One missing aspect of our understanding of $\Psi_{\tilde D}$ in (\ref{comm:PiO}) is whether or not
it is a homeomorphism.  By Proposition \ref{prop:cont}, this is equivalent to whether or not $\bP_{\tilde D}$ is 
surjective.  With the added hypothesis that $\bG$ admits a quasilogarithm and that $p^{r+1} > 2dim(\bG)$,
we prove this surjectivity of $\bP_{\tilde D}$ in Theorem \ref{thm:quasi-log}.
This enables us to conclude in Theorem \ref{thm:tD-test} that $M \mapsto \Pi(\tilde D(\Gr)_M$ detects
projectivity: a finite dimensional $\tilde D(\Gr)$-module $M$ is projective if and only if 
$\Pi(\tilde D(\Gr))_M$ is empty.
\vskip .1 in
 
We first recall the definition of a quasilogarithm for $\bG$.

\begin{defn} \cite{K-V}, \cite[Defn 6.2]{F-N}
\label{defn:quasi-log}
Let $\bG$ be a linear algebraic group over $k$.
A quasilogarithm for $\bG$ is a $\bG$-equivariant map of $k$-schemes
\begin{equation}
\label{eqn:quasi-log}
L: \bG \ \to \  \Spec S^\bu(\fg^*) \equiv  \ul{\fg}
\end{equation}
such that $L(e_\bG) = 0$ and such that the differential $d_{e_\bG}L: T_{e_\bG} \bG \to T_0 \ul{\fg} = \fg$ 
is the identity.

We denote by \ $I_r \ \subset \ S^\bu(\fg^*)$ the $\bG$-stable ideal generated by $\{ f^{p^r}, \ f \in \fg^* \}$.
\end{defn}

\vskip .1in

\begin{ex} \cite[Lem C3]{B-K-V}, \cite[6.4]{F-N}
\label{ex:quasi}
If $\bG$ is a simple algebraic group for which $p$ is very good, then $\bG$ admits a quasilogarithm.
Furthermore any Borel subgroup of $\bG$ also admits a quasilogarithm.  If $\bU$ is a unipotent 
subgroup of a semi-simple algebraic group $\bG$ which is normalized by a maximal torus and if 
$p$ is greater than the nilpotent class of $\bU$, then $\bU$ admits a quasilogarithm.
\end{ex}

\vskip .1in

For us, a critical consequence of the condition that $\bG$ admits a quasilogarithm
$L: \bG \to \ul \fg$  is the following result, a summary of the discussion prior to and 
including Lemma 6.6 of \cite{F-N}.

\begin{prop}
\label{prop:critical}
Assume that $\bG$ is a linear algebraic group smooth over $k$ which admits a quasilogarithm.
The $\bG$-equivariant map on coordinate algebras associated to the composition
 $L \circ i: \Gr \subset \bG \to \ul \fg$ 
has the form 
$$S^\bu(\fg^*) \ \to\ S^\bu(\fg^*)/I_r \ \stackrel{\sim}{\to} \  k[\Gr].$$ 

Moreover, the isomorphism $S^\bu(\fg^*)/I_r \ \stackrel{\sim}{\to} \  k[\Gr]$ 
induces an isomorphism of $k$-algebras
$$S^\bu(\fg^*)/I_r \# k\Gr \ \stackrel{\sim}{\to} \ k[\Gr]\# k\Gr \ \simeq  D(\Gr).$$

Similarly, $L$ determines the isomorphism of $k$-algebras
$$S^\bu(\fg^*)/I_{r+1} \# k\Gr \ \stackrel{\sim}{\to} \ k[\bG_{(r+1)}]\# k\Gr \ \simeq \tilde D(\Gr).$$

These smash products algebras $S^\bu(\fg^*)/I_r \# k\Gr, \  S^\bu(\fg^*)/I_{r+1} \# k\Gr$ are associated to the
coadjoint action of $\Gr$ on $\fg^*$.  They both admit a $\bG_{(r)}$-equivariant grading, with $\fg^*$ 
homogenous of degree 1 and $k\Gr$ homogeneous of degree 0. 
\end{prop}

\vskip .1in

Using the structure of a quasilogarithm on $\bG$, we obtain the following strengthening of Corollary \ref{cor:finite}.
The proof of this theorem is essentially that of \cite[Thm 6.9]{F-N} which gives the analogous 
result for $D(\bG_{(r)})$.

\begin{thm}
\label{thm:quasi-log}
Assume that $\bG$ is a linear algebraic group smooth over $k$ which admits a quasilogarithm.
Further assume that $p^{r+1} > 2dim(\bG)$.
Then the map 
$$\tilde \theta_r: H^*(\bG_{(r)},k) \otimes S^\bu(\fg^{(r+1)}[2]) \ \to H^*(\tilde D(\Gr),k)$$
of Theorem \ref{thm:fingen}
is a split injection with $H^*(\bG_{(r)},k) \otimes S^\bu(\fg^{(r+1)}[2])$-module complement 
consisting entirely of nilpotent elements. 

Similarly, the map
$i_O^*: H^*(\tilde D(\Gr),k) \ \to \ H^*(O(\Gr),k)$ has nilpotent kernel in $H^*(\tilde D(\Gr),k)$; and 
the image in $H^*(O(\Gr),k)$ has $H^*(\bG_{(r)},k) \otimes S^\bu(\fg^{(r+1)}[2])$-module complement 
consisting entirely of nilpotent elements.

Thus, $\bP_O:  \bP (H^*(O(\Gr),k)) \to \bP (H^*(\tilde D(\Gr),k)) $ is a  homeomorphism.  Consequently,
$\Psi_{\tilde D}: \Pi(\tilde D(\Gr)) \ \to \ \bP(H^\bu(\tilde D(\Gr),k))$ is also a homeomorphism
\end{thm}

\begin{proof}
The statement that $\tilde \theta_r$ is a split injection is that of Theorem \ref{thm:fingen}(1).

The gradings on $S^\bu(\fg^*)/I_r \# k\Gr, \  S^\bu(\fg^*)/I_{r+1} \# k\Gr$ (and thus also on 
$O(\Gr)$) determine internal gradings on the spectral sequences of Theorem \ref{thm:fingen}; we view
these gradings as $\bZ/p^{r+1}$-gradings.
As in \cite[Lem 6.7]{F-N}, the internal grading on $\fg^{(r+1)}[2] \subset  H^2(k[\bG_{(r+1)}],k)$ is $p^{r+1}$,
hence $\ \equiv 0 \ mod \ p^{r+1}$ and
the internal grading of $\fg[1] \subset H^1(k[\bG_{(r+1)}],k)$ is 1.   

We conclude that the map $\tilde \theta_r$ of Theorem \ref{thm:fingen} is an
isomorphism onto the degree $\equiv 0$ portion of the cohomoloy of $\tilde D(\Gr)$.  As argued in the proof of
\cite[Thm 6.9]{F-N}, the hypothesis that $p^{r+1} > 2dim(\bG)$ implies that this image has 
an $H^*(\bG_{(r)},k) \otimes S^\bu(\fg^{(r+1)}[2])$-module complement in 
$H^*(\tilde D(\Gr),k)$ which is contained in the nilradical.  

We may repeat this argument with the spectral sequence (\ref{eqn:Grothss}) for $\tilde D(\Gr)$ mapping to the
spectral sequence 
\begin{equation}
\label{eqn:GrothO}
E_2^{s,t} = H^s(\bG_{(r)}, H^t(k[G_{(r+1)}/\Gr],k)) \ \Rightarrow \ H^{s+t}(O(\bG_{(r)}),k).
\end{equation}  
Thus, the composition 
$H^*(\bG_{(r)},k) \otimes S^\bu(\fg^{(r+1)}[2]) \ \to \ H*(\tilde D(\Gr),k) \to \ H^*(O(\Gr),k)$
is also a split injection with $H^*(\bG_{(r)},k) \otimes S^\bu(\fg^{(r+1)}[2])$-module complement 
consisting entirely of nilpotent elements.  This readily implies the second statement 
concerning $i_O^*: H^*(\tilde D(\Gr),k) \ \to \ H^*(O(\Gr),k)$.  Consequently, $\bP_O$ is a 
bijection as well as a closed embedding, and thus is a homeomorphism.
This, together with Proposition \ref{prop:cont} implies that $\Psi_{\tilde D}$ is also a 
homeomorphism.
\end{proof}

We ``recall" the following projectivity criterion for a finite dimensional $\tilde D(\bG_{(r)})$-module
in terms of cohomological support varieties.

\begin{prop} (see \cite[Prop 1.5]{F-Par})
\label{prop:Hempty}
Let $M$ be a finitely generated $\tilde D(\bG_{(r)})$-module.  Then $M$ is projective if and only if
\ $\bP (H^\bu(\tilde D(\bG_{(r)}),k))_{ M}  \ = \ \emptyset.$

Similarly, if $M$ is a finitely generated $O(\Gr)$-module, them $M$ is projective if and only if 
\ $\bP (H^\bu(O(\bG_{(r)}),k))_{M}  \ = \ \emptyset.$
\end{prop}

\begin{proof}
As shown in \cite{L-S} for any finite dimensional Hopf algebra with finitely generated cohomology,
$\tilde D(\bG_{(r)})$ and 
$O(\Gr)$ are Frobenius algebras.
This enables the proof of \cite[Prop 1.5]{F-Par} to apply, even though that result was stated only
for Hopf algebras which are restricted enveloping algebras.
\end{proof}

Combining Theorem \ref{thm:quasi-log},  with Proposition \ref{prop:Hempty}, we conclude that
the support theory $M \mapsto  \Pi(\tilde D(\Gr))_M$ provides a test for projectivity.

\begin{thm}
\label{thm:tD-test}
Assume the hypotheses on $\bG$ of Theorem \ref{thm:quasi-log} and consider a finite dimensional
$\tilde D(\Gr)$-module $M$.  Then $M$ is projective as a $\tilde D(\Gr)$-module  if and only if  
$\Pi(\tilde D(\Gr))_M \ = \ \emptyset.$
\end{thm}

\begin{proof}
Assume that $M$ is projective as a $\tilde D(\Gr)$-module.   Then the pull-back along
$(\alpha_K + \beta_K)^*$ of $M$ is projective (and thus free) as a $K[t]/t^p$-module
for any $\pi$-point pair $(\alpha_K,\beta_K)$, so that $\Pi(\tilde D(\Gr))_M \ = \ \emptyset$.

Conversely, assume that  $\Pi(\tilde D(\Gr))_M \ = \ \emptyset$.
Since all the maps of the commutative square (\ref{comm:PiO}) are bijective (in fact,
homeomorphisms), we 
conclude that $\bP H^\bu((\tilde D(\Gr),k)_M \ = \ \emptyset.$  Thus, by Proposition 
\ref{prop:Hempty}, $M$ is a projective $\tilde D(\Gr)$-module.
\end{proof}

The arguments of Sections \ref{sec:cohom} and \ref{sec:quasi-log} apply to 
$\tilde D(G^\prime) \equiv k[\bG_{(r+1)}]\# G^\prime \ \subset \ \tilde D(\bG_{(r)})$
with merely a change of notation, as we state explicitly in the following corollary (and use in
the next section).

\begin{cor}
\label{cor:Gprime}
Let $G^\prime \ \hookrightarrow \ \Gr$ be a subgroup scheme of $\Gr$.  
Theorems \ref{thm:fingen} and \ref{thm:quasi-log} remain valid upon replacing $i_O: O(\Gr) \to \tilde D(\Gr)$
 by
  $$i_{O^\prime}:  k[\bG_{(r+1)}/\Gr] \# kG^\prime \  \hookrightarrow k[\bG_{(r+1)}] \# kG^\prime.$$
 
Consequently, as in Proposition \ref{prop:Hempty} and Theorem \ref{thm:tD-test}, a finite dimensional \\
$k[\bG_{(r+1)}] \# kG^\prime$-module $M$ is projective if and only if its restriction along $i_{O^\prime}$ is a projective
$ k[\bG_{(r+1)}/\Gr] \# kG^\prime$-module.
\end{cor}

\vskip .2in


\section{Relating $\Pi(\tilde D(\Gr))_M$ to $\bP(H^\bu(\tilde D(\Gr),k))_M$}
\label{sec:relate}

In this section, we verify for certain finite dimensional $\tilde D(\Gr)$-modules $M$ that the 
natural inclusion  $\Psi_{\tilde D,M}: \Pi(\tilde D(\Gr))_M \ \hookrightarrow \ \bP(H^\bu(\tilde D(\Gr),k))_M$
 is a bijection and thus a homeomorphism.
 
 We begin with an important class of examples (due to J. Carlson; see \cite{Ca}) of modules constructed 
to have given cohomological support.   Let  $P_* \to k$ be a minimal resolution of $k$ as a
$\tilde D(\Gr)$-module and define $\Omega_{\tilde D(\bG_{(r)})}^n(k)$ for $n \geq 1$ to be the kernel of 
$P_{n-1} \to P_{n-2}$ (where $P_{-1}$ is set equal to $k$).  
Recall that $\zeta \in H^n(\tilde D(\bG_{(r)}),k)$ naturally
corresponds to a map $\tilde \zeta: \Omega_{\tilde D(\bG_{(r)})}^n(k) \to k$ in the stable 
module category $stmod(\tilde D(\Gr))$.
The $\tilde D(\Gr)$-module  $L_\zeta$ is defined by the short exact sequence of $\tilde D(\bG_{(r)})$-modules
\begin{equation}
\label{eqn:Lzeta}
0 \to L_\zeta \ \to \ \Omega_{\tilde D(\bG_{(r)})}^n(k) \ \stackrel{\tilde \zeta}{\to} \  k \to 0.
\end{equation}
For $\xi \in H^n(O(\bG_{(r)}),k)$ with corresponding map $\tilde \zeta: \Omega_{O(\bG_{(r)})}^nk) \to k$ in the stable 
module category $stmod(O(\Gr))$, we similarly define the $O(\Gr)$-module $L_\xi$ by the 
short exact sequence of $O(\Gr)$-modules
$$0 \to L_\xi \ \to \ \Omega_{O\bG_{(r)})}^{n}(k) \ \stackrel{\tilde \xi}{\to} \  k \to 0.$$

We restrict attention to $n = 2d > 0$.   For such $n = 2d$, 
 $\Omega_{k[t]/t^p}^{2d}(k)$ is stably isomorphic to $k$.

\begin{prop}
\label{prop:Carlson}
Consider the $\tilde D(\Gr)$-module $L_\zeta$ associated to a cohomology class
$\zeta \in H^{2d}(\tilde D(\bG_{(r)}),k)$.  Then a  $\pi$-point pair $(\alpha_K,\beta_K)$ of $\tilde D(\bG_{(r)})$
satisfies the condition $(\alpha_K+\beta_K)^*(L_\zeta)$ is projective as a $K[t]/t^p$-module if and only if 
$(\alpha_K+\beta_K)^*(\zeta) \not= 0 \in H^{2d}(K[t]/t^p,K)$.

Consequently, $\Pi(\tilde D(\bG_{(r)}))_{L_\zeta}$ consists of equivalence classes of $\pi$-point pairs $(\alpha_K,\beta_K)$
such that $(\alpha_K+\beta_K)^*(\zeta) = 0 \in H^{2d}(K[t]/t^p,K)$.
\end{prop}

\begin{proof}
Carlson's identification of the support variety of $L_\zeta$ in the context of a group algebra of a finite group 
applies to $\tilde D(\Gr)$, for his argument merely 
requires that the algebra $\tilde D(\bG_{(r)})$ be a Frobenius algebra and that
$\alpha_K+\beta_K: K[t]/t^p \to \tilde D(\Gr)$ be flat.  The key observation of Carlson's proof is 
that the pullback along $\alpha_K+ \beta_K$ of (\ref{eqn:Lzeta}) is split if and only if
$(\alpha_K+\beta_K)^*(\zeta) \not= 0 \in H^{2d}(K[t]/t^p,K)$.
\end{proof}

Using Proposition \ref{prop:cont} and Theorem \ref{thm:quasi-log},  we verify that $\Psi_{\tilde D,L_\zeta}$
is a homeomorphism.

\begin{prop}
\label{prop:iso-Lz}
Assume the hypotheses on $\bG$ of Theorem \ref{thm:quasi-log}. 
Consider $\zeta \in H^\bu(\tilde D(\Gr),k)$ and denote by $\zeta_O$ the restriction
of $\zeta$ to  $H^\bu(O(\Gr),k)$.  Then
\begin{enumerate}
\item
The restriction of the $\tilde D(\Gr)$-module $L_\zeta$ to $O(\Gr)$ is stably equivalent to $L_{\zeta_O}$.
\item
$\bP i_O^*$ restricts to a bijection 
$$\bP_O: \bP(H^\bu(O(\Gr),k))_{L_{\zeta_O}} \ \stackrel{\sim}{\to} 
\ \bP(H^\bu(\tilde D(\Gr),k))_{L_{\zeta}}.$$
\item
$\Psi_{\tilde D}$ restricts to a homeomorphism
$$\Psi_{\tilde D,L_\zeta}: \Pi(\tilde D(\Gr))_{L_\zeta}  \ \stackrel{\sim}{\to} \ \bP(H^\bu(\tilde D(\Gr),k))_{L_\zeta}.$$
\end{enumerate}
\end{prop}

\begin{proof}
To prove the first assertion, observe that the $O(\Gr)$-module $\Omega_{O(\Gr}^{2d}(k)$ is stably equivalent to 
the restriction of $\Omega_{\tilde D(\Gr)}^{2d}(k)$ because the restriction of $\tilde D(\Gr)$-projectives to $O(\Gr)$
are projective $O(\Gr)$-modules.   Thus, restricting the short exact sequence 
(\ref{eqn:Lzeta}) of $\tilde D( \Gr)$-modules along $i_O: O(\Gr) \to \tilde D(\Gr)$ takes the form 
of a distinguished triangle
$$  L_{\zeta_O} \ \to \ \Omega_{O(\Gr}^{2d}(k)\  \stackrel{\zeta_O}{\to} \ k $$
in the stable module category $stmod(O(\Gr)$.

Using Proposition \ref{prop:Carlson}, we see that $\bP(H^\bu(O(\Gr),k))_{L_{\zeta_O}} \ \subset
\ \bP(H^\bu(O(\Gr),k))$ is the zero locus of the 
homogeneous ``function" $\zeta_O$ on $\bP(H^\bu(O(\Gr),k))$, and similarly that
$\bP(H^\bu(\tilde D(\Gr),k))_{L_{\zeta}} \ \subset \ \bP(H^\bu(\tilde D(\Gr),k))$ is the zero locus of the 
homogeneous ``function" $\zeta$ on $\bP(H^\bu(\tilde D(\Gr),k))$.  The second assertion thus follows from
the fact that $\zeta_O \ = \ i_O^*(\zeta)$ and the fact that $\bP_O$ is a homeomorphism by Theorem
\ref{thm:quasi-log}.  

The last assertion now follows from the commutative square (\ref{comm:PiO}) of Proposition \ref{prop:cont} 
and the fact that $\Psi_O$ restricts to a bijection $\Psi_{O,L_{\zeta_O}}: \Psi(\tilde O(\Gr))_{L_{\zeta_O}} \ 
\stackrel{\sim}{\to}\ \bP(H^\bu(\tilde O(\Gr),k))_{L_{\zeta_O}}$ by Proposition \ref{prop:PsiO}.
\end{proof}

The following proposition gives some insight into the topology on $\Pi(\tilde D(\Gr))$ beyond 
knowing that it is the topology inherited from $\bP (H^\bu(\tilde D(\bG_{(r)}),k))$.

\begin{prop}
\label{prop:bases}
Assume the hypotheses on $\bG$ of Theorem \ref{thm:quasi-log} so that the maps of the
commutative square (\ref{comm:PiO}) are homeomorphisms.  Then 
$$\{ \Pi(\tilde D(\Gr))_{L_\zeta}, \ \zeta \in H^\bu(\tilde D(\Gr),k) \ {\text homogeneous} \}$$
is a closed base for the topological space $\Pi(\tilde D(\Gr))$.
\end{prop}

\begin{proof}
Observe that
$$\{ \bP(H^\bu(O(\Gr),k))_{L_\xi} \ \xi \in H^\bu(O(\Gr),k) \ {\text homogeneous} \}$$
is a closed base for the topology of $\bP(H^\bu(O(\Gr),k))$.  Thus, the fact that
$\bP_O: \bP(H^\bu(O(\Gr),k)) \ \to \ \bP(H^\bu(\tilde D(\Gr),k))$ is a homeomorphism 
(by Theorem \ref{thm:quasi-log}) together with  
Proposition \ref{prop:iso-Lz}(2) implies that 
$$\{ \bP(H^\bu(\tilde D(\Gr),k))_{L_\zeta}, \ \zeta \in H^\bu(D(\Gr),k) \ {\text homogeneous} \}$$
is a closed base for the topology of $\bP(H^\bu(\tilde D(\Gr),k))$.  Since $\Psi_{\tilde D}$ is a 
homeomorphism, Proposition \ref{prop:iso-Lz}(3) implies that 
$$\{ \Pi(\tilde D(\Gr))_{L_\zeta}, \ \zeta \in H^\bu(\tilde D(\Gr),k) \ {\text homogeneous} \}$$
is a closed base for the topological space $\Pi(\tilde D(\Gr))$.
\end{proof}

	We next consider the class of  $\tilde D(\Gr)$-modules given as 
the inflation along $\tilde D(\Gr) \to \Gr$ of finite dimensional $\Gr$-modules
(equivalently, as finite dimensional $\tilde D(\Gr)$-module whose restriction to 
$k[\bG_{(r+1)}] \subset \tilde D(\Gr)$ is trivial).

\begin{prop}
\label{prop:restrict-triv}
Assume the hypotheses on $\bG$ of Theorem \ref{thm:quasi-log}.  Let $M$ be a finite dimensional
 $\tilde D(\Gr)$-module whose restriction to $k[\bG_{(r+1)}] \subset \tilde D(\Gr)$ is trivial.  Then
\begin{enumerate}
\item  $\bP_{O,M}: \bP(H^\bu(O(\Gr),k))_M \ \to \ \bP(H^\bu(\tilde D(\Gr),k))_M$ is a homeomorphism.
\item
$\bP(H^\bu(O(\Gr),k))_M$ can be identified as the ``geometric join" of $\Proj(S^\bu(\fg^{(r+1)}))$ and $\bP (H^\bu(\Gr,k))_M$.
\item
 $\Psi_{\tilde D,M}: \Pi(\tilde D(\Gr))_M \ \to \ \bP(H^\bu(\tilde D(\Gr),k))_M$ is also a homeomorphism.
\end{enumerate}
\end{prop}

\begin{proof}
We consider the commutative square
\begin{equation}
\label{diag:annihil}
\begin{xy}*!C\xybox{%
\xymatrix{
H^\bu(\tilde D(\Gr),k)  \ar[r]^{i_O^*} \ar[d] & H^\bu(O(\Gr),k)  \ar[d] \\
H^\bu(\tilde D(\Gr),M^*\otimes M) \ar[r]^{i_O^*}  & H^\bu(O(\Gr),M^*\otimes M). }
}\end{xy}
\end{equation}

We observe that there are multiplicative spectral sequences for $M^*\otimes M$ replacing $k$ corresponding to
the spectral sequences (\ref{eqn:Grothss}) and (\ref{eqn:GrothO}).  Equipping $M^*\otimes M$ with 
internal grading degree 0, (\ref{diag:annihil}) determines a commutative square of multiplicative spectral sequences 
with internal grading whose $E_2$-page has the form
\begin{equation}
\label{diag:E2pages}
\begin{xy}*!C\xybox{%
\xymatrix{
H^s(\Gr,H^t(k[\bG_{(r+1)}],k))  \ar[r] \ar[d] & H^s(\Gr,H^t(k[\bG_{(r+1)}/\Gr],k)) \ar[d] \\
H^s(\Gr,H^t(k[\bG_{(r+1)}],k) \otimes M^* \otimes M)  \ar[r] & 
H^s(\Gr,H^t(k[\bG_{(r+1)}/\Gr],k) \otimes M^* \otimes M) .}
}\end{xy}
\end{equation}

The restriction map $H^*(k[\bG_{(r+1)}],k) \to H^*(k[\bG_{(r+1)}/\Gr],k)$ has nilpotent kernel 
and cokernel which are nilpotent, $\Gr$-summands of $H^*(k[\bG_{(r+1)}/\Gr],k)$.  Thus, arguing as 
in the proof of Theorem \ref{thm:quasi-log}, classes at the $E_2$-level with internal grading 
not congruent to 0 modulo $p^{r+1}$ are 
torsion; moreover, the lower horizontal arrow of (\ref{diag:E2pages}) restricted to internal degrees 
congruent to 0 modulo $p^{r+1}$ has nilpotent kernel 
and cokernel which are nilpotent, $\Gr$-summands.

Since $\bP_O$ is a homeomorphism, to prove (1) it suffices to verify that if $\fp \subset H^\bu(O(\Gr),k)$ 
is a homogeneous prime ideal whose inverse image $\fq \equiv i_O^{-1}(\fp)$ contains the 
kernel of the left vertical map of (\ref{diag:annihil}), then $\fp$ contains the kernel of the 
right vertical map.   The validity of this statement is equivalent to the same statement of the commutative
square obtained from (\ref{diag:annihil}) by replacing the lower horizontal arrow
by its associated map on reduced algebras
(i.e., by dividing out be the nilradicals of the domain and range of this map).  By the preceding 
discussion, this latter map is an isomorphism.  Thus, having made
this replacement, the required statement about $i_O^{-1}(\fp)$ is immediate.

Assertion (2) follows from the familiar interpretation of the projectivization of the product of two affine varieties.

Assertion (3) follows from the commutativity of (\ref{comm:PiO}) and assertion (1).
\end{proof}

We consider a third class of finite dimensional $\tilde D(\Gr)$-modules, those
with the property that their restrictions to $k[\bG_{(r+1)}]$ are projective.  We first recall the 
following detection theorem of A. Suslin whose formulation we quote from \cite[Thm 4.10]{F-P1}.

\begin{thm} \cite{Suslin}
\label{thm:suslin}
Let $G$ be a finite group scheme, $\Lambda$ a unital associative $G$-algebra, and $\zeta \in H^\bu(G,k)$
be a homogeneous cohomology class of even degree.  Then $\zeta$ is nilpotent if and only if $\zeta_K$
restricts to a nilpotent class in the cohomology of every quasi-elementary subgroup scheme $\cE_K \subset G_K$
for any field extension $K/k$.
\end{thm}

Let $M$ be a finite dimensional $\tilde D(\Gr)$-module with the property that the restriction of
$M$ to $k[\bG_{(r+1)}]$ is projective.   If $M$ satisfies this property, then the spectral sequences 
(\ref{eqn:Grothss}) and (\ref{eqn:GrothM} degenerate so that the edge homomorphisms provide the
 isomorphisms of algebras
\begin{equation}
\label{eqn:GrH0}
H^*(\Gr,H^0(k[\bG_{(r+1)}],M^*\otimes M)) \ \stackrel{\sim}{\to} \ H^*(\tilde D(\Gr),M^*\otimes M),
\end{equation}
$$H^*(\Gr,H^0(k[\bG_{(r+1)}/\Gr],M^*\otimes M)) \ \stackrel{\sim}{\to} \ H^*(O(\Gr),M^*\otimes M).$$
We use (\ref{eqn:GrH0}) in conjunction with Theorem \ref{thm:suslin} to prove that 
$$\Psi_{\tilde D_M}: \Pi(\tilde D(\Gr))_M \ \to \ \bP (H^\bu(\tilde D(\Gr),k))_M$$
is a homeomorphism for these modules.

\begin{prop}
\label{prop:tD-proj}
Assume the hypotheses on $\bG$ of Theorem \ref{thm:quasi-log}.
Let $M$ be a finite dimensional $\tilde D(\Gr)$-module with the property that 
the restriction of $M$ to $k[\bG_{(r+1)}]$ is projective.
Then the restriction of $\Psi_{\tilde D}$ is a homeomorphism
$$\Psi_{\tilde D,M}: \Pi(\tilde D(\Gr))_M \  \stackrel{\sim}{\to} \ \bP (H^\bu(\tilde D(\Gr),k))_M.$$.
\end{prop}

\begin{proof}
Let $\Lambda(M)$ denote unital associative $\Gr$-algebra $H^0(k[\bG_{(r+1)}],M^*\otimes M)$ and let 
$\Lambda_O(M)$ denote $H^0(k[\bG_{(r+1)}/\Gr],M^*\otimes M) $.  Using 
the isomorphisms (\ref{eqn:GrH0}), we rewrite (\ref{diag:annihil}) as
\begin{equation}
\label{diag:Suslinsquare}
\begin{xy}*!C\xybox{%
\xymatrix{
H^\bu(\tilde D(\Gr),k)  \ar[r]^{i_O^*} \ar[d] & H^\bu(O(\Gr),k)  \ar[d] \\
H^\bu(\Gr,\Lambda(M)) \ar[r]^{i_O^*}  & H^\bu(\Gr,\Lambda_O(M)). }
}\end{xy}
\end{equation}
We set $I_M$ to be the kernel of the left vertical map and $J_M$ the kernel of the right vertical map 
of (\ref{diag:Suslinsquare}).

Since we know that $\Psi_{\tilde D,M}$ is well defined and injective and that $\Psi_{\tilde D}$ is bijective,
to prove the proposition we must show for any homogeneous prime ideal $\fq \subset H^\bu(O(\Gr),k)$ with 
inverse image $\fp \subset H^\bu(\tilde D(\Gr),k)$ that $J_M \subset \fq$ (i.e., $\fp$ is a point of 
$\bP(H^\bu(O(\Gr),k))_M$) whenever $I_M \subset \fp$.  Assume that $\zeta \in \fp \subset H^\bu(\tilde D(\Gr),k)$
is a homogeneous element which is not in $I_M$, so that its image $\ol \zeta \in H^\bu(\Gr,\Lambda(M))$ is not nilpotent.  
Using Theorem \ref{thm:suslin}, let $\gamma: \bG_{a(r),K} \to \bG_{(r),K}$ be a 1-parameter subgroup with
the property that $\gamma^*(\ol \zeta) \in H^\bu(\bG_{a(r),K},\gamma^*(\Lambda(M_K)))$ is not  nilpotent. 
  
In other words, $\gamma^*(\zeta) \in H^\bu(K[\bG_{(r+1)}]\# \bG_{a(r)},k)$ has image 
$\ol{\gamma^*(\zeta)} \in Ext_{K[\bG_{(r+1)}]\# \bG_{a(r),K}]}^*(\gamma^*(M_K),\gamma^*(M_K))$ which is not nilpotent. 
By Corollary \ref{cor:Gprime}, the image of 
$\ol{\gamma^*(\zeta)}$ in $Ext_{K[\bG_{(r+1)}/\Gr]\# \bG_{a(r),K}]}^*(\gamma^*(M_K),\gamma^*(M_K))$
 is not nilpotent as well.  Consequently, $i_O^*(\zeta)$ does not lie in $J_M$.  
 Thus, we conclude if $\fp$ contains $I_M$ then $\fq$ must contain $J_M$.

\end{proof}

\vskip .2in


\section{Further remarks}
\label{sec:remarks}

We begin this final section with a natural question.

\begin{question}
Assume that $\bG$ satisfies the conditions of Theorem \ref{thm:quasi-log}.
Can one extend the results of Section \ref{sec:relate} to prove that
$$\Psi_{\tilde D(\bG_{(r)}),M}: \Pi(\tilde D(\bG_{(r)}))_M \ \to \
\bP(H^\bu(\tilde D(\bG_{(r)}),k)_M$$ 
is a homeomorphism for all finite dimensional $\tilde D(\bG_{(r)})$-modules?
\end{question}

We proceed to  briefly introduce Hopf subalgebras $\tilde D^{(s)}(\Gr) \ \subset \tilde D(\Gr)$
(generalizing $O(\Gr) \subset \tilde D(\Gr)$)
and quotient maps of Hopf algebras $k[\bG_{r+s}] \ \twoheadrightarrow \ D(\Gr)$
(generalizing $\tilde D(\Gr) \twoheadrightarrow D(\Gr))$.   These generalizations behave much as the algebras
we have consider in previous sections.

\begin{prop}
\label{prop:twistedD}
The  $\bG$-equivariant quotient map 
$\bG_{(r+1)}\ \twoheadrightarrow\  \bG_{(r+1)}/\bG_{(1)} \simeq (\bG^{(1)})_{(r)}$ determines a Hopf algebra
embedding of smash products
\begin{equation}
\label{eqn:iD}
i_D: \tilde D^{(1)}(\bG_{(r)})\ \equiv k[(\bG^{(1)})_{(r)}] \# k\Gr  \ \hookrightarrow\  
k[\bG_{(r+1}] \# k\Gr \equiv \tilde D(\Gr).
\end{equation}
This is a map of free (right) $O(\Gr)$-modules

The composition 
$$\tilde D^{(1)}(\bG_{(r)})  \ \hookrightarrow \  \tilde D(\Gr) \ \ \twoheadrightarrow \ D(\Gr)$$
is the map $F \# id: k[(\bG^{(1)})_{(r)}]\# k\Gr \ \to \ k[\bG_{(r)}] \# k\Gr$, where $F$ is the 
Frobenius map.
\end{prop}

\begin{proof}
As in the justification of (\ref{eqn:extDrinfeld}),  the map of (\ref{eqn:iD}) is an embedding of
Hopf algebras.  We choose a $k$-linear section $\sigma: k[\Gr/\bG_{(1)}] \to k[\bG_{(r+1)}/\bG_{(1)}]$
of the quotient map induced by the embedding $\Gr/\bG_{(1)} \subset \bG_{(r+1)}/\bG_{(1)}$ and 
consider the 
$$k[(\bG_{(r+1)}/\bG_{(1)})/(\Gr/\bG_{(1)})] \ \simeq \ k[\bG_{(r+1)}/\bG_{(r)}]\ \simeq \ k[(\bG^{(r)})_{(1)}]$$
bilinear map
$$ k[(\bG^{(1)})_{(r-1)}] \otimes  k[(\bG^{(r)})_{(1)}]  \ \to \ k[(\bG^{(1)})_{(r)}]$$
analogous to (\ref{eqn:sig}).  Since $k[(\bG^{(1)})_{(r-1)})]\otimes 1$ is central in $\tilde D^{(1)}(\bG_{(r)})$,
we conclude as in the proof of Proposition \ref{prop:Oflat} that this provides $\tilde D^{(1)}(\bG_{(r)})$
with the structure of a free (right) $O(\Gr)$-module.

To verify that the asserted composition equals $F \# id$, we use the identification 
$\bG_{(r+1)}/\bG_{(1)} \simeq (\Gr)^{(1)}$.
\end{proof}

We can view $\tilde D^{(1)}(\bG_{(r)})$
as a twisted form of $D(\Gr)$ for which the action of $\Gr$ on $k[(\bG^{(1)})_{(r)}]$ equals the Frobenius
twist of the coadjoint action of $\Gr$ on $k[\Gr]$.

The proof of the following theorem is a basically identical to the proofs given in previous sections of
the corresponding properties for $\tilde D(\Gr)$-modules.

\begin{prop}
\label{prop:tildeD1}
For any $\tilde D^{(1)}(\bG_{(r)})$-module $M$, we can define $M \mapsto \Pi(\tilde D^{(1)}(\bG_{(r)}))_M$
in strict analogy with $M \mapsto \Pi(\tilde D(\bG_{(r)}))_M$.  So defined, $\pi$-point pairs of 
$\tilde D^{(1)}(\bG_{(r)})$ are $\pi$-point pairs of $O(\Gr)$ with the equivalence relation determined by
consideration of finite dimensional $\tilde D^{(1)}(\bG_{(r)})$-modules.  For any finite dimensional
$\tilde D^{(1)}(\bG_{(r)})$-module $M$, we define $\Pi(\tilde D^{(1)}(\bG_{(r)}))_M \subset \Pi(\tilde D^{(1)}(\bG_{(r)}))$
to consist of  those equivalence classes of $\pi$-point pairs $(\alpha_K,\eta_K)$ with the property that 
$(\alpha_K + \beta_K)^*(M_K)$ is not projective.

There is a natural continuous map  $\Psi_{\tilde D^{(1)}(\bG_{(r)}}: \Pi(\tilde D^{(1)}(\bG_{(r)})) \to
\bP(H^\bu(\tilde D^{(1)}(\bG_{(r)}),k)$.
Assuming that $\bG$ admits a quasilogarithm and $p^{r+1} > 2dim(\bG)$, this map is a 
homeomorphism and restricts to a homeomorphism 
$$\Psi_{\tilde D^{(1)}(\bG_{(r)},M}: \Pi(\tilde D^{(1)}(\bG_{(r)}))_M \ \stackrel{\sim}{\to} \
\bP(H^\bu(\tilde D^{(1)}(\bG_{(r)}),k)_M$$
provided that $M$ is in one of the classes analogous to those classes of 
 $\tilde D(\bG_{(r)})$-modules considered in Section \ref{sec:relate}.
\end{prop}

\vskip .1in

\begin{remark}
The triple of Hopf algebras over $k$
$$O(\Gr) \ \hookrightarrow \ \tilde D^{(1)}(\Gr)  \ \hookrightarrow \ \tilde D(\Gr)$$
can be refined to a sequence of embeddings of Hopf algebras which are free over $O(\Gr)$
$$O(\Gr) \equiv D^{(r)}(\Gr) \ \hookrightarrow \ \cdots
 \hookrightarrow \tilde D^{(1)}(\Gr)  \hookrightarrow \tilde D^{(0)}(\Gr) \equiv \tilde D(\Gr),$$
where
$$\tilde D^{(s)}(\bG_{(r)}) \ \equiv \ k[(\bG^{(s)}_{(r+1-s)}] \# k\Gr, \quad 0 \leq s \leq r.$$
One readily checks that replacing $\tilde D^{(1)}(\bG_{(r)})$ in Proposition \ref{prop:tildeD1} 
by $\tilde D^{(s)}(\bG_{(r)})$ gives analogous properties for $\tilde D^{(s)}(\bG_{(r)})$-modules.
\end{remark}

\vskip .1in

\begin{remark}
For any $s \geq 1$, we define the Hopf subalgebra
$$k[\bG_{(r+s)}] \# k\Gr \ \subset \ k[\bG_{(r+s)}] \# k\bG_{(r+s)} \equiv D(\bG_{(r+s)}.$$
This Hopf algebra admits a quotient map of Hopf algebras
$$k[\bG_{(r+s)}] \# k\Gr \ \twoheadrightarrow \ D(\bG_{(r))}$$
determined by the quotient map $k[\bG_{(r+s)}] \twoheadrightarrow k[\bG_{(r+s)}] \to k[\Gr]$
as in (\ref{eqn:quot}).  We obtain a sequence of quotient Hopf algebras
$$ \cdots \twoheadrightarrow k[\bG_{(r+s)}] \# k\Gr \twoheadrightarrow k[\bG_{(r+s-1)}] \# k\Gr 
\twoheadrightarrow \cdots \to \tilde D(\Gr) \twoheadrightarrow D(\Gr).$$

We leave to the reader to check that our results for $\tilde D(\bG_{(r)})$-modules extend
to $k[\bG_{(r+s)}] \# k\Gr$ for $s \geq 1$.
\end{remark}

\vskip .1in
 \begin{remark}
We recall that a theorem of Pevtsova, Suslin, and the author \cite[Thm 4.10] {FPS} asserts that 
if $M$ is a finite dimensional $G$-module for a finite group scheme $G$ over $k$
 and if $\alpha_K: K[t]/t^p \to G_K$ is a $\pi$-point of $G$ at which the Jordan type
 $\alpha_K^*(M_K)$ is maximal among all Jordan types $\beta_L^*(M_L)$ 
 as $\beta_L: L[t]/^p \to G_L$ varies among all $\pi$-points of $G$, then the Jordan type
 of $\alpha_K^{\prime *}(M_K)$ equals the Jordan type of $\alpha_K^*(M_K)$ 
 whenever the $\pi$-point $\alpha_K^\prime$ is equivalent to $\alpha_K$.
 
 This immediately applies to $O(\Gr)$-modules, since $O(\Gr)$ is isomorphic as a $k$-algebra
 to the group algebra of a finite group scheme as seen in Proposition \ref{prop:ulg}.
 By the definition of the equivalence relation on $\pi$-point pairs of $\tilde D(\Gr)$,
 this remains valid for a $\tilde D(\Gr)$-module $M$ whose maximal Jordan type
 is that of a projective module.  
  \end{remark}
 
 \begin{question}
 Does this ``independence of representative
 of equivalence class" remain valid for any finite dimensional $\tilde D(\Gr)$-module $M$?  Namely,
 if the  equivalence class of a $\pi$-point pair of $\tilde D(\Gr)$ contains some
 some representative $(\alpha_K,\beta_K)$ at which $(\alpha_K + \beta_K)^*(M_K)$ has 
 maximal Jordan type, then does $(\alpha_L^\prime + \beta_L)^*(M_L)$ have this
 same Jordan type for any $(\alpha_K^\prime + \beta_K^\prime) \ 
 \sim \ (\alpha_K + \beta_K)$?
\end{question}

\vskip .2in


\end{document}